# EFFICIENT ESTIMATION OF COPULA-BASED SEMIPARAMETRIC MARKOV MODELS

By Xiaohong Chen[1], Wei Biao Wu[2] and Yanping Yi

*Yale University, University of Chicago and New York University*

*This paper is dedicated to Professor Peter C. B. Phillips on the occasion of his 60th birthday*

This paper considers the efficient estimation of copula-based semiparametric strictly stationary Markov models. These models are characterized by nonparametric invariant (one-dimensional marginal) distributions and parametric bivariate copula functions where the copulas capture temporal dependence and tail dependence of the processes. The Markov processes generated via tail dependent copulas may look highly persistent and are useful for financial and economic applications. We first show that Markov processes generated via Clayton, Gumbel and Student's $t$ copulas and their survival copulas are all geometrically ergodic. We then propose a sieve maximum likelihood estimation (MLE) for the copula parameter, the invariant distribution and the conditional quantiles. We show that the sieve MLEs of any smooth functional is root-$n$ consistent, asymptotically normal and efficient and that their sieve likelihood ratio statistics are asymptotically chi-square distributed. Monte Carlo studies indicate that, even for Markov models generated via tail dependent copulas and fat-tailed marginals, our sieve MLEs perform very well.

**1. Introduction.** A copula function is a multivariate probability distribution function with uniform marginals. A copula-based method has become one popular tool for modeling nonlinearity, asymmetricality and tail dependence in financial and insurance risk managements. See Embrechts, McNeil and Straumann (2002), McNeil, Frey and Embrechts (2005), Embrechts

Received December 2008; revised March 2009.
[1]Supported by the NSF Grant SES-0631613.
[2]Supported by the NSF Grant DMS-04-78704.

([2009](#)), Genest, Gendron and Bourdeau-Brien ([2008](#)), Patton ([2002](#), [2006](#), [2008](#)) and the references therein for reviews of various theoretical properties and financial applications of the copula approach.

While the majority of the previous work using copulas has focused on modeling the contemporaneous dependence between multiple univariate series, there are also a growing number of papers using copulas to model the temporal dependence of a univariate nonlinear time series. Granger ([2003](#)) defines persistence (such as "long memory" or "short memory") for general nonlinear time series models via copulas. Darsow, Nguyen and Olsen ([1992](#)), de la Peña, Ibragimov and Sharakhmetov ([2006](#)) and Ibragimov ([2009](#)) provide characterizations of a copula-based time series to be a Markov process. Joe ([1997](#)) proposes a class of parametric (strictly) stationary Markov models based on parametric copulas and parametric invariant (one-dimensional marginal) distributions. Chen and Fan ([2006](#)) study a class of semiparametric stationary Markov models based on parametric copulas and nonparametric invariant distributions.

Let $\{Y_t\}$ be a stationary Markov process of order one with a continuous invariant (one-dimensional marginal) distribution $G$. Then its probabilistic properties are completely determined by the bivariate joint distribution function of $Y_{t-1}$ and $Y_t$, $H(y_1, y_2)$ (say). By Sklar's theorem [see McNeil, Frey and Embrechts ([2005](#)), Nelsen ([2006](#))], one can uniquely express $H(\cdot,\cdot)$ in terms of the invariant distribution $G$ and the bivariate copula function $C(\cdot,\cdot)$ of $Y_{t-1}$ and $Y_t$,

$$H(y_1, y_2) \equiv C(G(y_1), G(y_2)).$$

Thus one can always specify a stationary first-order Markov model with continuous state space by directly specifying the marginal distribution of $Y_t$ and the bivariate copula function of $Y_{t-1}$ and $Y_t$. The advantage of the copula approach is that one can freely choose the marginal distribution and the bivariate copula function separately; the former characterizes the marginal behavior such as the fat-tails and/or skewness of the time series $\{Y_t\}_{t=1}^n$ while the latter characterizes all the temporal dependence properties that are invariant to any increasing transformations as well as the tail dependence properties of the time series. Although being strictly stationary first-order Markov, a model generated via a copula (especially a tail-dependent copula) is very flexible. This model can generate a rich array of nonlinear time series patterns, including persistent clustering of extreme values via tail dependent copulas evaluated at fat-tailed marginals, asymmetric dependence, and other "look alike" behaviors present in many popular nonlinear models such as ARCH, GARCH, stochastic volatility, near-unit root, long-memory, models with structural breaks, Markov switching and so on. From the point of view of financial applications, one attractive property of the copula-based Markov model is that the implied conditional quantiles are automatically monotonic across quantiles. This nice feature has been exploited by Chen, Koenker and Xiao (2008) and Bouyé and Salmon (2008) in their study of copula-based nonlinear quantile autoregression and value at risk (VaR).

In this paper, we shall focus on the class of copula-based, strictly stationary, semiparametric first-order Markov models, in which the true copula density function has a parametric form ($c(\cdot, \cdot; \alpha_0)$), and the true invariant distribution is of an unknown form ($G_0(\cdot)$) but is absolutely continuous with respect to the Lebesgue measure on the real line. Any model of this class is completely described by two unknown characteristics: the copula dependence parameter $\alpha_0$ and the invariant distribution $G_0(\cdot)$. To establish the asymptotic properties of any semiparametric estimators of $(\alpha_0, G_0)$, one needs to know temporal dependence properties of the copula-based Markov models. For this class of models, Chen and Fan (2006) show that the $\beta$-mixing temporal dependence measure is purely determined by the properties of the copulas (and does not depend on the invariant distributions); and Beare (2008) provides simple sufficient conditions for geometric $\beta$-mixing in terms of copulas without any tail dependence [such as Gaussian, Frank and Eyraud–Farlie–Gumbel–Morgenstern (EFGM) copulas]. Neither paper is able to verify whether or not a Markov process generated via a tail dependent copula (such as Clayton, survival Clayton, Gumbel, survival Gumbel or Student's $t$) is geometric $\beta$-mixing. Ibragimov and Lentzas (2008) demonstrate via simulation that Clayton copula-based first-order strictly stationary Markov models could behave as "long memory" in copula levels. In this paper, we show that Clayton, survival Clayton, Gumbel, survival Gumbel and Student's $t$ copula-based Markov models are actually geometrically ergodic (hence geometric $\beta$-mixing). Therefore, according to our theorem, although a time series plot of a Clayton copula (or survival Clayton, Gumbel, survival Gumbel or other tail-dependent copula) generated Markov model may look highly persistent and "long memory alike," it is, in fact, weakly dependent and "short memory."

In this paper, we propose a sieve maximum likelihood estimation (MLE) procedure for the copula parameter $\alpha_0$, the invariant distribution $G_0$ and the conditional quantiles of a copula-based semiparametric Markov model. This procedure approximates the unknown marginal density by flexible parametric families of densities with increasing complexity (sieves), and then maximizes the joint likelihood with respect to the unknown copula parameter and the sieve parameters of the approximating marginal density. We show that the sieve MLEs of any smooth functionals of $(\alpha_0, G_0)$ are root-$n$ consistent, asymptotically normal and efficient; and that their sieve likelihood ratio statistics are asymptotically chi-square distributed. We also present simple consistent estimators of asymptotic variances of the sieve MLEs of smooth functionals. It is interesting to note that although the conditional distribution of a copula-based semiparametric stationary Markov model depends on the unknown invariant distribution, the plug-in sieve MLE estimators of the nonlinear conditional quantiles (VaR) are still $\sqrt{n}$-consistent, asymptotically normal and efficient.

To the best of our knowledge, Atlason (2008) is the only other paper that also considers the semiparametric efficient estimation of a copula parameter $\alpha_0$ for a copula-based first-order strictly stationary Markov model. His work and ours were done at the same time, but independently. While we propose the sieve likelihood joint estimation of $G_0$ and $\alpha_0$, Atlason (2008) proposes the rank likelihood estimation of the copula parameter $\alpha_0$, and relies on a simulation method to evaluate his rank likelihood. However, Atlason (2008) does not investigate the semiparametric efficient estimation of the invariant distribution $G_0$ nor the conditional quantiles.

Previously, Chen and Fan (2006) proposed a simple two-step estimation procedure in which one first estimates the invariant CDF $G_0(\cdot)$ by a re-scaled empirical CDF $G_n$ of the data $\{Y_t\}_{t=1}^n$, and then estimates the copula parameter $\alpha_0$ by maximizing the pseudo log-likelihood corresponding to copula density evaluated at pseudo observations $\{G_n(Y_t)\}_{t=1}^n$. Chen and Fan's procedure can be viewed as an extension of the one proposed by Genest, Ghoudi and Rivest (1995) for a bivariate copula-based joint distribution model of a random sample $\{(X_i, Y_i)\}_{i=1}^n$ to a univariate first-order Markov model of a time series data $\{Y_i\}_{i=1}^n$ (with $X_i = Y_{i-1}$). Both are semiparametric analogs of the two-step parametric procedure that is called the "inference functions for margins" (IFM) in Joe (1997), Chapter 10. Just as the two-step estimator of Genest, Ghoudi and Rivest (1995) is generally inefficient for a bivariate random sample [see, e.g., Genest and Werker (2002)], the two-step estimator of Chen and Fan (2006) is inefficient for a univariate Markov model.

We present Monte Carlo studies to compare the finite sample performance of our sieve MLE, the two-step estimator of Chen and Fan (2006), the correctly specified parametric MLE and the incorrectly specified parametric MLE for Clayton, Gumbel, Gaussian, Frank and EFGM copula-based Markov models. Numerous simulation studies demonstrate that the two-step estimator of Chen and Fan (2006) is not only inefficient but also severely biased (in finite sample) when the time series has strong tail dependence, and it leads to a biased and inefficient plug-in estimator of conditional quantiles (or VaR). The simulation results indicate that our sieve MLEs perform very well; when the copula-based Markov process has strong tail dependence, the sieve MLEs have much smaller biases and smaller variances than the two-step estimators.

The rest of this paper is organized as follows. In Section 2, we present the class of copula-based semiparametric strictly stationary Markov models and show that many widely used tail dependent copula-based Markov models are geometrically $\beta$-mixing. In Section 3, we introduce the sieve MLE, and obtain its consistency and rate of convergence. Section 4 establishes the asymptotic normality and semiparametric efficiency of the sieve MLE. Section 5 shows that their sieve likelihood ratio statistics are asymptotically chi-square distributed which suggests a simple way to construct confidence regions for the copula parameter and other smooth functionals. In Section 6, we first briefly review some popular existing estimators. We then conduct some simulation studies to compare the finite sample performance of our sieve MLE and these alternative estimators. Section 7 briefly concludes. All the proofs are relegated to the Appendix.

Finally, we wish to point out that given the characterization results of Darsow, Nguyen and Olsen (1992) and Ibragimov (2009) on higher order Markov models via copulas, we can easily extend our sieve MLE method and our results for copula-based first-order Markov models to copula-based higher order Markov models. For presentational clarity we do not give the details here.

**2. Copula-based Markov models.** In this section, we first present the model and then some implied temporal dependence properties.

2.1. *The model.* Darsow, Nguyen and Olsen (1992) provide characterization of first-order Markov processes by bivariate copulas and one-dimensional marginal distributions (see Nelsen [(2006), Section 6.4] for a brief review). Throughout this paper, we assume that the true data generating process (DGP) satisfies the following assumption:

ASSUMPTION M. (DGP): (1) $\{Y_t : t = 1, \ldots, n\}$ is a sample of a strictly stationary first-order Markov process generated from $(G_0(\cdot), C(\cdot, \cdot; \alpha_0))$ where $G_0(\cdot)$ is the true invariant distribution that is absolutely continuous with respect to the Lebesgue measure on the real line (with its support $\mathcal{Y}$, a nonempty interval of $\mathcal{R}$); $C(\cdot, \cdot; \alpha_0)$ is the true parametric copula for $(Y_{t-1}, Y_t)$ up to unknown value $\alpha_0$ and is absolutely continuous with respect to the Lebesgue measure on $[0,1]^2$. (2) The true marginal density $g_0(\cdot)$ of $G_0(\cdot)$ is positive on its support $\mathcal{Y}$; and the true copula density $c(\cdot, \cdot; \alpha_0)$ of $C(\cdot, \cdot; \alpha_0)$ is positive on $(0,1)^2$.

In Assumption M(1), the assumption of absolute continuity of the bivariate copula $C(\cdot, \cdot; \alpha_0)$ rules out the Fréchet–Hoeffding upper ($C(u_1, u_2) = \min(u_1, u_2)$) and the lower ($C(u_1, u_2) = \max(u_1 + u_2 - 1, 0)$) bounds, as well as their linear combinations [and, say, shuffles and Min copulas discussed in Darsow, Nguyen and Olsen (1992)].

Under Assumption M(1), the true conditional probability density function, $p^0(\cdot|Y^{t-1})$ of $Y_t$ given $Y^{t-1} \equiv (Y_{t-1}, \ldots, Y_1)$, is given by

$$(2.1) \qquad p^0(\cdot|Y^{t-1}) = h_0(\cdot|Y_{t-1}) \equiv g_0(\cdot)c(G_0(Y_{t-1}), G_0(\cdot); \alpha_0),$$

where $h_0(\cdot|Y_{t-1})$ denotes the true conditional density of $Y_t$, given $Y_{t-1}$. Under Assumption M(1), the transformed process $\{U_t : U_t \equiv G_0(Y_t)\}_{t=1}^n$ is also a strictly stationary first-order Markov process with uniform marginals and $C(\cdot, \cdot; \alpha_0)$, the joint distribution of $U_{t-1}$ and $U_t$. Then $C_{2|1}[\cdot|u; \alpha_0] \equiv \frac{\partial}{\partial u} C(u, \cdot; \alpha_0) \equiv C_1(u, \cdot; \alpha_0)$ is the conditional distribution of $U_t \equiv G_0(Y_t)$, given $U_{t-1} = u$; and $C_{2|1}^{-1}[q|u; \alpha_0]$ is the $q$th, $q \in (0,1)$, conditional quantile of $U_t$, given $U_{t-1} = u$.

Note that the conditional density of $Y_t$, given $Y^{t-1}$, is a function of both the copula density $c(\cdot, \cdot; \alpha_0)$ and the marginal density $g_0$; hence the $q$th, $q \in (0,1)$, the conditional quantile of $Y_t$ given $Y^{t-1}$ is also a function of both the copula and the marginal

$$(2.2) \qquad Q_q^Y(y) = G_0^{-1}(C_{2|1}^{-1}[q|G_0(y); \alpha_0]).$$

By definition, $C_{2|1}^{-1}[q|u; \alpha_0]$ is increasing in $q$; hence the $q$th conditional quantile of $Y_t$ given $Y^{t-1}$, $Q_q^Y(y)$, is also increasing in $q$.

2.2. *Tail dependence, temporal dependence.* All the dependence measures that are invariant under increasing transformations can be expressed in terms of copulas [see, e.g., McNeil, Frey and Embrechts (2005) and Nelsen (2006)]. For example, Kendall's tau is $\tau = 4 \iint_{[0,1]^2} C(u_1, u_2) \, dC(u_1, u_2) - 1$, and Spearman's rho is $\rho_S = 12 \iint_{[0,1]^2} C(u_1, u_2) \, du_1 \, du_2 - 3$. The lower (resp. upper) tail dependence coefficients $\lambda_L$ (resp. $\lambda_U$) in terms of copulas are

$$\lambda_L \equiv \lim_{u \to 0^+} \Pr(U_2 \leq u | U_1 \leq u) = \lim_{u \to 0^+} \frac{C(u,u)}{u} \quad \text{and}$$

$$\lambda_U \equiv \lim_{u \to 1^-} \Pr(U_2 \geq u | U_1 \geq u) = \lim_{u \to 1^-} \frac{1 - 2u + C(u,u)}{1 - u},$$

provided the limits exist. [See Kortschak and Albrecher (2009) for examples of copulas with nonexisting limits for tail dependence and their applications.]

For financial risk management, the Markov models generated via tail-dependent copulas are much more relevant than models without tail dependence. In particular, the following three examples have been widely used in financial applications:

EXAMPLE 2.1 (Clayton copula-based Markov model). The bivariate Clayton copula is

$$C(u_1, u_2, \alpha) = [u_1^{-\alpha} + u_2^{-\alpha} - 1]^{-1/\alpha}, \qquad 0 \leq \alpha < \infty.$$

Clayton copula has Kendall's tau $\tau = \frac{\alpha}{2+\alpha}$, and the lower tail dependence coefficient $\lambda_L = 2^{-1/\alpha}$ that is increasing in $\alpha$, but no upper tail dependence. Clayton copula becomes the independence copula $C_I(u_1, u_2) = u_1 u_2$ in the limit when $\alpha \to 0$.

EXAMPLE 2.2 (Gumbel copula-based Markov model). The bivariate Gumbel copula is

$$C(u_1, u_2; \alpha) = \exp(-[(-\ln u_1)^\alpha + (-\ln u_2)^\alpha]^{1/\alpha}), \qquad 1 \leq \alpha < \infty.$$

Gumbel copula has Kendall's tau $\tau = 1 - \frac{1}{\alpha}$, and the upper tail dependence coefficient $\lambda_U = 2 - 2^{1/\alpha}$ that is increasing in $\alpha$, but no lower tail dependence. Gumbel copula becomes the independence copula $C_I(u_1, u_2) = u_1 u_2$ in the limit when $\alpha \to 1$.

EXAMPLE 2.3 (Student $t$ copula-based Markov model). The bivariate Student $t$ copula is

$$C(u_1, u_2; \alpha) = \mathbf{t}_{\nu,\rho}(t_\nu^{-1}(u_1), t_\nu^{-1}(u_2)), \qquad \alpha = (\nu, \rho), |\rho| < 1, \nu \in (1, \infty],$$

where $\mathbf{t}_{\nu,\rho}(\cdot, \cdot)$ is the bivariate Student-$t$ distribution with mean zero, the correlation matrix having off-diagonal element $\rho$, and degrees of freedom $\nu$, and $t_\nu(\cdot)$ is the CDF of a univariate Student-$t$ distribution with mean zero, and degrees of freedom $\nu$. Student $t$ copula has Kendall's tau $\tau = \frac{2}{\pi} \arcsin \rho$, and symmetric tail dependence $\lambda_L = \lambda_U = 2t_{\nu+1}(-\sqrt{(\nu+1)(1-\rho)/(1+\rho)})$ that is decreasing in $\nu$. The Student $t$ copula becomes a Gaussian copula in the limit when $\nu \to \infty$.

2.2.1. *Geometric $\beta$-mixing.* For analyzing asymptotic properties of any semiparametric estimators of $(\alpha_0, G_0)$, it is convenient to apply empirical process results for strictly stationary geometrically ergodic (or geometric $\beta$-mixing) Markov processes. See Appendix A for some equivalent definitions of $\beta$-mixing and ergodicity for strictly stationary Markov processes.

REMARK 2.1. (1) Under Assumption M, the time series $\{Y_t\}_{t=1}^n$ is strictly stationary ergodic and is also $\beta$-mixing (see, e.g., Bradley [(2005), Corollary 3.6] and Chen and Fan (2006)).

(2) Proposition 2.1 of Chen and Fan (2006) presents high-level sufficient (and almost necessary) conditions in terms of a copula to ensure $\beta$-mixing decaying either exponentially fast or polynomially fast. Their working-paper version points out that their Proposition 2.1 implies the Markov models based on Gaussian and EFGM copulas are geometric $\beta$-mixing. However, they do not verify whether any other copulas satisfy the conditions of their Proposition 2.1.

(3) Beare [(2008), Theorem 3.1 and Remark 3.5] shows that all Markov models generated via symmetric absolute continuous copulas with positive and square integrable copula densities are geometric $\beta$-mixing. In Remark 3.7, he points out that many commonly used bivariate copulas without tail dependence, such as Gaussian, EFGM, Frank, Gamma, binomial and hypergeometric copulas, satisfy the conditions of his Theorem 3.1.

(4) Beare [(2008), Theorem 3.2] shows that all bivariate absolute continuous copulas with square integrable densities do not have any tail dependence. Although he shows that a Markov model based on Student's $t$ copula is rho mixing and hence is geometrically strong mixing, Beare (2008) does not verify whether a Markov model generated via any tail dependent copula (such as Clayton, Gumbel or Student's $t$ copula) is geometrically $\beta$-mixing.

Ibragimov and Lentzas (2008) demonstrate via simulation that Clayton copula generated first-order strictly stationary Markov models behave as "long memory" processes in copula levels when the Clayton copula parameter $\alpha$ is big. The time series plots (see Figure 1) of such Markov processes appear to be "long memory alike." (See Section 6.2 on how to simulate copula-based first-order stationary Markov time series. The clusterings of extremes in Figure 1 are due to tail dependence properties of Clayton and Gumbel copulas.) Nevertheless, our next theorem shows that they are in fact geometrically ergodic and hence they are "short memory" processes.

THEOREM 2.1 (Geometric ergodicity). *Under Assumption M, the Markov time series $\{Y_t\}_{t=1}^n$ generated via Clayton copula with $0 < \alpha < \infty$, Gumbel copula with $1 \leq \alpha < \infty$, Student's $t$ copula with $|\rho| < 1$ and $2 \leq \nu < \infty$, are all geometrically ergodic (and hence geometrically $\beta$-mixing).*

REMARK 2.2. If $\{U_t\}_{t=1}^n$ is a $C_U(\cdot,\cdot)$ copula generated strictly stationary first-order Markov model with uniform marginals, then $\{V_t \equiv 1 - U_t\}_{t=1}^n$ is also a copula-based strictly stationary first-order Markov model with uniform marginals and bivariate copula function

$$C_V(v_1, v_2) \equiv \Pr(V_{t-1} \leq v_1, V_t \leq v_2) = \Pr(U_{t-1} \geq 1 - v_1, U_t \geq 1 - v_2)$$
$$= v_1 + v_2 - 1 + C_U(1 - v_1, 1 - v_2) \equiv C_U^s(v_1, v_2)$$

which is the survival copula of $C_U^s(u_1, u_2)$ [see Nelsen (2006)]. Therefore, a copula $C_U(\cdot, \cdot)$ generated strictly stationary first-order Markov process is geometrically ergodic or $\beta$-mixing with certain decay speed $\beta_j = o(1)$ if and only if its survival copula $C_U^s(\cdot, \cdot)$ generated Markov process is geometrically ergodic or $\beta$-mixing with the same decay speed $\beta_j = o(1)$.

By Theorem 2.1 and Remark 2.2, we immediately see that survival Clayton and survival Gumbel generated first-order stationary Markov processes are also geometrically ergodic.

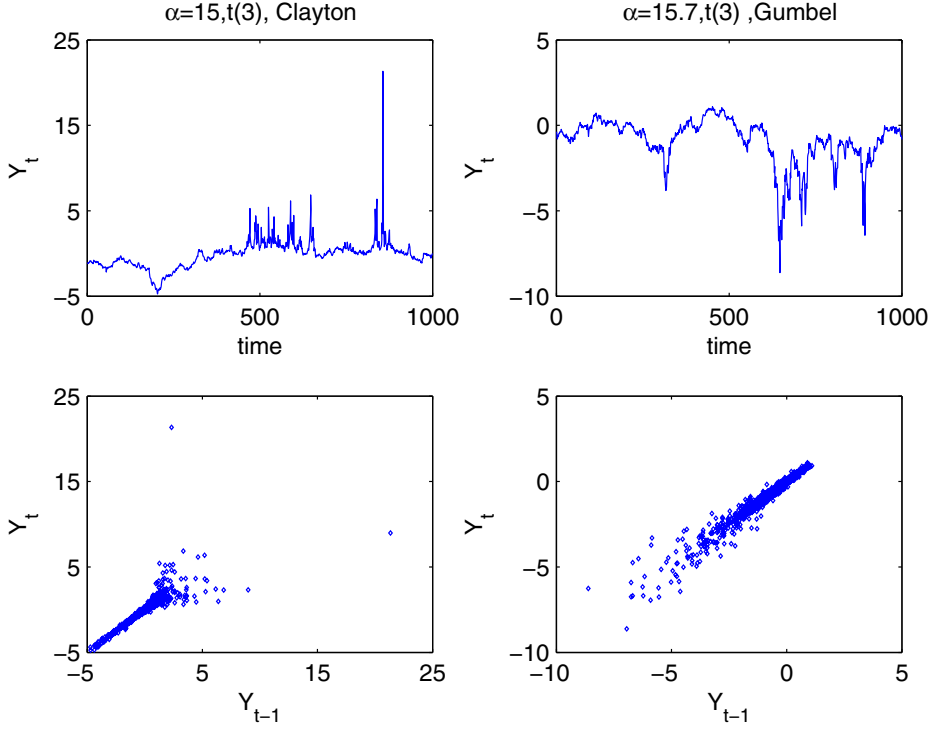

Fig. 1. *Markov time series: tail dependence index* $= 0.9548$, *Student* $t_3$ *marginal distribution.*

**3. Sieve MLE, consistency with rate.** Under Assumption M, we see that the true conditional density $p^0(\cdot|Y^{t-1})$ of $Y_t$ given $Y^{t-1} \equiv (Y_{t-1}, \ldots, Y_1)$ is given by (2.1). Let

$$p(\cdot|Y^{t-1}) = h(\cdot|Y_{t-1}; \alpha, g) \equiv g(\cdot)c(G(Y_{t-1}), G(\cdot); \alpha)$$

denote any candidate conditional density of $Y_t$ given $Y^{t-1}$. Let $Z_t = (Y_{t-1}, Y_t)$, and denote

$$\begin{aligned}\ell(\alpha, g, Z_t) &\equiv \log p(Y_t|Y^{t-1}) = \log\{h(Y_t|Y_{t-1}; \alpha, g)\} \\ &\equiv \log g(Y_t) + \log c(G(Y_{t-1}), G(Y_t); \alpha) \\ &\equiv \log g(Y_t) + \log c\bigg(\int 1(y \leq Y_{t-1})g(y)\,dy, \int 1(y \leq Y_t)g(y)\,dy; \alpha\bigg)\end{aligned}$$

as the log-likelihood associated with the conditional density $p(Y_t|Y^{t-1})$. Here $1(\cdot)$ stands for the indicator function. Then the joint log-likelihood function of the data $\{Y_t\}_{t=1}^n$ is given by

$$L_n(\alpha, g) \equiv \frac{1}{n}\sum_{t=2}^n \ell(\alpha, g, Z_t) + \frac{1}{n}\log g(Y_1).$$

The approximate sieve MLE $\widehat{\gamma}_n \equiv (\widehat{\alpha}_n, \widehat{g}_n)$ is defined as

$$(3.1) \qquad L_n(\widehat{\alpha}_n, \widehat{g}_n) \geq \max_{\alpha \in \mathcal{A}, g \in \mathcal{G}_n} L_n(\alpha, g) - O_p(\delta_n^2),$$

where $\delta_n$ is a positive sequence such that $\delta_n = o(1)$, and $\mathcal{G}_n$ denotes the sieve space [i.e., a sequence of finite dimensional parameter spaces that become dense (as $n \to \infty$) in the entire parameter space $\mathcal{G}$ for $g_0$].

There exist many sieves for approximating a univariate probability density function. In this paper, we will focus on using linear sieves to directly approximate either a square root density:

$$(3.2) \qquad \mathcal{G}_n = \left\{ g_{K_n} \in \mathcal{G} : g_{K_n}(y) = \left[\sum_{k=1}^{K_n} a_k A_k(y)\right]^2, \int g_{K_n}(y)\,dy = 1 \right\},$$

$$K_n \to \infty, \frac{K_n}{n} \to 0;$$

or a log density:

$$(3.3) \qquad \mathcal{G}_n = \left\{ g_{K_n} \in \mathcal{G} : g_{K_n}(y) = \exp\left\{\sum_{k=1}^{K_n} a_k A_k(y)\right\}, \int g_{K_n}(y)\,dy = 1 \right\},$$

$$K_n \to \infty, \frac{K_n}{n} \to 0,$$

where $\{A_k(\cdot) : k \geq 1\}$ consists of known basis functions, and $\{a_k : k \geq 1\}$ is the collection of unknown sieve coefficients.

Suppose the support $\mathcal{Y}$ (of the true $g_0$) is either a compact interval (say $[0, 1]$) or the whole real line $\mathcal{R}$. Let $r > 0$ be a real-valued number, and $[r] \geq 0$ be the largest integer such that $[r] < r$. A real-valued function $g$ on $\mathcal{Y}$ is said to be $r$-smooth if it is $[r]$ times continuously differentiable on $\mathcal{Y}$, and its $[r]$th derivative satisfies a Hölder condition with exponent $r - [r] \in (0, 1]$ (i.e., there is a positive number $K$ such that $|D^{[r]}g(y) - D^{[r]}g(y')| \leq K|y - y'|^{r-[r]}$ for all $y, y' \in \mathcal{Y}$. Here $D^{[r]}$ stands for the differential operator). We denote $\Lambda^r(\mathcal{Y})$ as the class of all real-valued functions on $\mathcal{Y}$ which are $r$-smooth; it is called a Hölder space.

Let the true marginal density function $g_0$ satisfy either $\sqrt{g_0} \in \Lambda^r(\mathcal{Y})$ or $\log g_0 \in \Lambda^r(\mathcal{Y})$. Then any function in $\Lambda^r(\mathcal{Y})$ can be approximated by some appropriate sieve spaces. For example, if $\mathcal{Y}$ is a bounded interval and $r > 1/2$, it can be approximated by the spline sieve $Spl(s, K_n)$ with $s > [r]$, the polynomial sieve, the trigonometric sieve, the cosine series and so on. When the support of $\mathcal{Y}$ is unbounded, thin-tailed density can be approximated by a Hermite polynomial sieve, while a polynomial fat-tailed density can be approximated by a spline wavelet sieve. See Chen (2007) for detailed descriptions of various sieve spaces $\mathcal{G}_n$. In our simulation study, we choose the sieve number in terms of $K_n$ using a modified AIC, although one could also use cross-validation [see, e.g., Fan and Yao (2003), Gao (2007), Li and Racine (2007)] and other computationally more intensive model selection methods [see, e.g., Shen, Huang and Ye (2004)] to choose the sieve number in terms of $K_n$. See Chen, Fan and Tsyrennikov (2006) for further discussions.

3.1. *Consistency.* In the following, we denote $Q_n(\alpha, g) \equiv \frac{n-1}{n} E_0[\ell(\alpha, g, Z_2)] + \frac{1}{n} E_0[\log g(Y_1)]$ where $E_0$ is the expectation under the true DGP (i.e., Assumption M). Denote $\gamma \equiv (\alpha, g)$ and $\gamma_0 \equiv (\alpha_0, g_0) \in \Gamma \equiv \mathcal{A} \times \mathcal{G}$.

ASSUMPTION 3.1. (1) $\alpha_0 \in \mathcal{A}$, where $\mathcal{A}$ is a compact set of $\mathcal{R}^d$ with a nonempty interior, $c(u_1, u_2; \alpha) > 0$ for all $(u_1, u_2) \in (0,1)^2$, $\alpha \in \mathcal{A}$; (2) $g_0 \in \mathcal{G}$, either $\mathcal{G} = \{g = f^2 > 0 : f \in \Lambda^r(\mathcal{Y}), \int g(y)\,dy = 1\}$ and $\mathcal{G}_n$ given in (3.2), or $\mathcal{G} = \{g = \exp(f) > 0 : f \in \Lambda^r(\mathcal{Y}), \int g(y)\,dy = 1\}$ and $\mathcal{G}_n$ given in (3.3), $r > 1/2$; (3) $Q_n(\alpha_0, g_0) > -\infty$, there is the metric $\|\gamma\|_c \equiv \sqrt{\alpha'\alpha} + \|g\|_c$ on $\Gamma \equiv \mathcal{A} \times \mathcal{G}$ and a positive measurable function $\eta(\cdot)$ such that for all $\varepsilon > 0$ and for all $k \geq 1$,

$$Q_n(\alpha_0, g_0) - \sup_{\alpha \in \mathcal{A}, g \in \mathcal{G}_k : \|\gamma_0 - \gamma\|_c \geq \varepsilon} Q_n(\alpha, g) \geq \eta(\varepsilon) > 0;$$

(4) the sieve spaces $\mathcal{G}_n$ are compact under the metric $\|g\|_c$; (5) there is $\Pi_n \gamma_0 \in \Gamma_n \equiv \mathcal{A} \times \mathcal{G}_n$ such that $\|\Pi_n \gamma_0 - \gamma_0\|_c = o(1)$; and $|Q_n(\Pi_n \gamma_0) - Q_n(\gamma_0)| = o(1)$.

For the norm $\|\gamma\|_c \equiv \sqrt{\alpha'\alpha} + \|g\|_c$ on $\Gamma \equiv \mathcal{A} \times \mathcal{G}$, one can use either a sup norm $\|g\|_\infty$ (or a weighted sup norm) or even a lower order Hölder norm $\|g\|_{\Lambda^{r'}}$ for $r' \in [0, r)$ (or its weighted version).

ASSUMPTION 3.2. (1) $E_0[\sup_{\gamma \in \Gamma_n} |\ell(\gamma, Z_t)|]$ is bounded; (2) there is a finite constant $\kappa > 0$ and a measurable function $M(\cdot)$ with $E_0[M(Z_t)] \leq const. < \infty$, such that for all $\delta > 0$,

$$\sup_{\{\gamma, \gamma_1 \in \Gamma_n : \|\gamma - \gamma_1\|_c \leq \delta\}} |\ell(\gamma, Z_t) - \ell(\gamma_1, Z_t)| \leq \delta^\kappa M(Z_t) \qquad \text{a.s.-}Z_t.$$

We note that under Assumption 3.1(1), (4), Assumption 3.2(1) is implied by Assumption 3.2(2).

PROPOSITION 3.1. *Under Assumptions M, 3.1 and 3.2, $\delta_n = o(1)$, $K_n \to \infty$ and $\frac{K_n}{n} \to 0$, we have*

$$\|\widehat{\gamma}_n - \gamma_0\|_c = o_p(1).$$

3.2. *Convergence rate.* Given the consistency result Proposition 3.1, $\varphi_n := \inf\{h > 0 \colon \Pr(\|\widehat{\gamma}_n - \gamma_0\|_c > h) \leq h\}$, the Levy distance between $\|\widehat{\gamma}_n - \gamma_0\|_c$ and 0 converges to 0. Let $\mathcal{N} = \{\gamma \in \Gamma \colon \|\gamma - \gamma_0\|_c \leq \varphi_n\}$ be the new parameter space, and the corresponding shrinking neighborhood in the sieve space, denoted as $\mathcal{N}_n = \mathcal{N} \cap \Gamma_n$, be the new sieve parameter space. Denote $\mathrm{Var}_0$ as the variance under the true DGP (i.e., Assumption M).

ASSUMPTION 3.3. (1) There are metric $\|\gamma\|_s \equiv \sqrt{\alpha'\alpha} + \|g\|_s$ on $\mathcal{N}$ such that $\|\gamma\|_s \leq \|\gamma\|_c$, and a constant $J_0 > 0$ such that for all $\varepsilon > 0$ and for all $n \geq 1$,

$$Q_n(\alpha_0, g_0) - \sup_{\gamma \in \mathcal{N}_n \colon \|\gamma_0 - \gamma\|_s \geq \varepsilon} Q_n(\alpha, g) \geq J_0 \varepsilon^2 > 0.$$

(2) $\sup_{\{\gamma \in \mathcal{N}_n \colon \|\gamma_0 - \gamma\|_s \leq \epsilon\}} \mathrm{Var}_0(\ell(\gamma, Z_t) - \ell(\gamma_0, Z_t)) \leq \mathit{const.} \times \epsilon^2$ for all small $\epsilon > 0$.

Assumption 3.3 suggests that a natural choice of $\|\gamma\|_s$ could be $(Q_n(\gamma_0) - Q_n(\gamma))^{1/2}$.

ASSUMPTION 3.4. (1) $\{Y_t\}_{t=1}^n$ is geometrically ergodic (hence geometrically $\beta$-mixing); (2) there is a constant $\kappa \in (0, 2)$ and a measurable function $M(\cdot)$ with $E_0[M(Z_t)^2 \log(1 + M(Z_t))] \leq \mathit{const.} < \infty$, such that for any $\delta > 0$,

$$\sup_{\{\gamma \in \mathcal{N}_n \colon \|\gamma_0 - \gamma\|_s \leq \delta\}} |\ell(\gamma, Z_t) - \ell(\gamma_0, Z_t)| \leq \delta^\kappa M(Z_t) \qquad \text{a.s.-}Z_t.$$

Although we do not need any $\beta$-mixing decay rates to establish consistency in Proposition 3.1, we need some $\beta$-mixing decay rates for rate of convergence.[3] Given the results in Section 2.2.1, Assumption 3.4(1) is typically satisfied by copula-based Markov models. Note that in Assumption 3.4(2), the moment restriction on the envelop function $M(Z_t)$ is weaker than the one $(E_0[M(Z_t)^\varsigma] \leq \mathit{const.} < \infty$ for some $\varsigma > 2)$ imposed in Chen and Shen (1998). This is because Chen and Shen (1998) only assume $\beta$-mixing with polynomial decay speed while our Assumption 3.4(1) assumes geometric $\beta$-mixing. It is well known that there are trade-offs between speed of mixing decay rate and finiteness of moments [see, e.g., Doukhan, Massart and Rio (1995) and Nze and Doukhan (2004)]. Assumption 3.4(2) is a very weak regularity condition and is satisfied whenever $\sup_{\eta \in [0,1], \gamma \in \mathcal{N}_n \colon \|\gamma_0 - \gamma\|_s \leq \delta} |\frac{d\ell(\gamma_0 + \eta[\gamma - \gamma_0], Z_t)}{d\eta}| \leq \delta^\kappa M(Z_t)$ with $M(Z_t)$ having a finite slightly higher than a second moment,

---

[3]It is common to assume some $\beta$-mixing or strong mixing decay rates in semi/nonparametric estimation and testing [see, e.g., Robinson (1983), Andrews (1994), Fan and Yao (2003), Gao (2007), Li and Racine (2007)].

which is satisfied by all the copula-based Markov models that satisfy the regularity conditions in Chen and Fan (2006) for semiparametric two-step estimators.

The next proposition is a direct application of Theorem 1 of Chen and Shen (1998), hence we omit its proof.

PROPOSITION 3.2. *Under Assumptions M, 3.1–3.4, we have*

$$\|\widehat{\gamma}_n - \gamma_0\|_s = O_p(\delta_n), \qquad \delta_n = \max\left\{\sqrt{\frac{K_n}{n}}, \|\gamma_0 - \Pi_n \gamma_0\|_s\right\} = o(1).$$

**4. Normality and efficiency of sieve MLE of smooth functionals.** Let $\rho: \mathcal{A} \times \mathcal{G} \to \mathcal{R}$ be a smooth functional and $\rho(\widehat{\gamma}_n)$ be the plug-in sieve MLE of $\rho(\gamma_0)$. In this section, we extend the results of Chen, Fan and Tsyrennikov (2006) on root-$n$ normality and efficiency of their sieve MLE for copula-based multivariate joint distribution model using i.i.d. data to our scalar strictly stationary first-order Markov setting.

4.1. $\sqrt{n}$-*asymptotic normality of* $\rho(\widehat{\gamma}_n)$. Recall that $\delta_n$ is the speed of convergence of $\|\widehat{\gamma}_n - \gamma_0\|_s$ to zero in probability, let $\mathcal{N}_0 = \{\gamma \in \mathcal{N} : \|\gamma_0 - \gamma\|_s \leq \delta_n \log \delta_n^{-1}\}$ and $\mathcal{N}_{0n} = \{\gamma \in \mathcal{N}_n : \|\gamma_0 - \gamma\|_s \leq \delta_n \log \delta_n^{-1}\}$, then $\widehat{\gamma}_n \in \mathcal{N}_{0n}$ with probability approaching one. Also denote $(U_1, U_2) = (G_0(Y_1), G_0(Y_2))$, $u = (u_1, u_2) \in [0, 1]^2$ and $c(G_0(Y_{t-1}), G_0(Y_t); \alpha_0) = c(U; \alpha_0) = c(\gamma_0, Z_t)$ (with the danger of slightly abusing notation).

ASSUMPTION 4.1. $\alpha_0 \in \text{int}(\mathcal{A})$.

ASSUMPTION 4.2. The second-order partial derivatives $\frac{\partial^2 \log c(u;\alpha)}{\partial \alpha \alpha'}$, $\frac{\partial^2 \log c(u;\alpha)}{\partial u_j \partial \alpha}$, $\frac{\partial^2 \log c(u;\alpha)}{\partial u_j \partial u_k}$ for $k, j = 1, 2$, are all well defined and continuous in $\gamma \in \mathcal{N}_0$.

Denote **V** as the linear span of $\Gamma - \{\gamma_0\}$. Under Assumption 4.2, for any $v = (v_\alpha, v_g)' \in \mathbf{V}$, we see that $\ell(\gamma_0 + \eta v, Z)$ is continuously differentiable in $\eta \in [0, 1]$. For any $\gamma \in \mathcal{N}_0$, define the first-order directional derivative of $\ell(\gamma, Z_t)$ at the direction $v \in \mathbf{V}$ as

$$\begin{aligned}
\frac{\partial \ell(\gamma, Z_t)}{\partial \gamma'}[v] &\equiv \frac{d\ell(\gamma + \eta v, Z_t)}{d\eta}\bigg|_{\eta=0} \\
&= \frac{\partial \log c(\gamma, Z_t)}{\partial \alpha'}[v_\alpha] + \frac{v_g(Y_t)}{g(Y_t)} \\
&\quad + \sum_{j=1}^{2} \frac{\partial \log c(\gamma, Z_t)}{\partial u_j} \int \mathbf{1}\{y \leq Y_{t-2+j}\} v_g(y) \, dy,
\end{aligned}$$

and the second-order directional derivative as

$$\frac{\partial^2 \ell(\gamma, Z_t)}{\partial \gamma \, \partial \gamma'}[v, \widetilde{v}] \equiv \frac{d}{d\widetilde{\eta}} \left\{ \frac{\partial \ell(\gamma + \widetilde{\eta}\widetilde{v}, Z_t)}{\partial \gamma'}[v] \right\} \bigg|_{\widetilde{\eta}=0}$$
$$= \frac{d^2 \ell(\gamma + \eta v + \widetilde{\eta}\widetilde{v}, Z_t)}{d\widetilde{\eta} \, d\eta} \bigg|_{\eta=0} \bigg|_{\widetilde{\eta}=0}.$$

ASSUMPTION 4.3. (1) $0 < E_0[(\frac{\partial \ell(\gamma_0, Z_t)}{\partial \gamma'}[v])^2] < \infty$ for $v \neq 0, v \in \mathbf{V}$;
(2) $\int \sup_{\eta \in \mathcal{S}_v} |\frac{dh(y|Y_{t-1};\gamma_0 + \eta v)}{d\eta}| \, dy < \infty$ and $\int \sup_{\eta \in \mathcal{S}_v} |\frac{d^2 h(y|Y_{t-1};\gamma_0 + \eta v)}{d\eta^2}| \, dy < \infty$ almost surely, for $\mathcal{S}_v = \{\eta \in [0,1] : \gamma_0 + \eta v \in \mathcal{N}_0\}$, $v \neq 0, v \in \mathbf{V}$.

Assumption 4.3(2) is a condition that is assumed even for parametric Markov models similar to those in Joe [(1997), Chapter 10] and Billingsley (1961b).

LEMMA 4.1.  *Under Assumptions M, 3.1(1), (2), 4.1, 4.2 and 4.3, we have, for any $v \in \mathbf{V}$, (1) $E_0((\frac{\partial \ell(\gamma_0, Z_t)}{\partial \gamma'}[v])(\frac{\partial \ell(\gamma_0, Z_s)}{\partial \gamma'}[\tilde{v}])) = 0$ for $\tilde{v} \in \mathbf{V}$ and all $s < t$. (2) $\{\frac{\partial \ell(\gamma_0, Z_t)}{\partial \gamma'}[v]\}_{t=1}^n$ is a martingale difference sequence with respect to the filtration $\mathcal{F}_{t-1} = \sigma(Y_1; \ldots; Y_{t-1})$. (3) $E_0((\frac{\partial \ell(\gamma_0, Z_t)}{\partial \gamma'}[v])^2) = -E_0(\frac{\partial^2 \ell(\gamma_0, Z_t)}{\partial \gamma \partial \gamma'}[v, v])$.*

Lemma 4.1 suggests that we can define the Fisher inner product on the space $\mathbf{V}$ as

$$\langle v, \tilde{v} \rangle \equiv E_0 \left[ \left( \frac{\partial \ell(\gamma_0, Z_t)}{\partial \gamma'}[v] \right) \left( \frac{\partial \ell(\gamma_0, Z_t)}{\partial \gamma'}[\tilde{v}] \right) \right]$$

and the Fisher norm for $v \in \mathbf{V}$ as $\|v\|^2 \equiv \langle v, v \rangle$. Let $\overline{\mathbf{V}}$ be the closed linear span of $\mathbf{V}$ under the Fisher norm. Then $(\overline{\mathbf{V}}, \|\cdot\|)$ is a Hilbert space.

The asymptotic properties of $\rho(\widehat{\gamma}_n)$ depend on the smoothness of the functional $\rho$ and the rate of convergence of $\widehat{\gamma}_n$. For any $v \in \overline{\mathbf{V}}$, we denote

$$\frac{d\rho(\gamma_0 + \eta v)}{d\eta} \bigg|_{\eta=0} \equiv \frac{\partial \rho(\gamma_0)}{\partial \gamma'}[v],$$

whenever the limit is well defined.

ASSUMPTION 4.4. (1) For any $v \in \overline{\mathbf{V}}$, $\rho(\gamma_0 + \eta v)$ is continuously differentiable in $\eta \in [0, 1]$ near $\eta = 0$, and

$$\left\| \frac{\partial \rho(\gamma_0)}{\partial \gamma'} \right\| \equiv \sup_{v \in \overline{\mathbf{V}}: \|v\| > 0} \frac{|\frac{\partial \rho(\gamma_0)}{\partial \gamma'}[v]|}{\|v\|} < \infty;$$

(2) there exist constants $c > 0$, $\omega > 0$, and a small $\epsilon > 0$ such that

$$\left| \rho(\gamma_0 + v) - \rho(\gamma_0) - \frac{\partial \rho(\gamma_0)}{\partial \gamma'}[v] \right| \leq c\|v\|^\omega \qquad \text{for any } v \in \overline{\mathbf{V}} \text{ with } \|v\| < \epsilon.$$

Under this assumption, by the Riesz representation theorem, there exists a $v^* \in \overline{\mathbf{V}}$ such that

(4.1) $$\frac{\partial \rho(\gamma_0)}{\partial \gamma'}[v] \equiv \langle v^*, v \rangle \qquad \text{for all } v \in \overline{\mathbf{V}}$$

and

$$\|v^*\|^2 = \left\| \frac{\partial \rho(\gamma_0)}{\partial \gamma'} \right\|^2 = \sup_{v \in \overline{\mathbf{V}}: \|v\| > 0} \frac{|\frac{\partial \rho(\gamma_0)}{\partial \gamma'}[v]|^2}{\|v\|^2} < \infty.$$

ASSUMPTION 4.5. (1) $\|\widehat{\gamma}_n - \gamma_0\| = O_p(\delta_n)$ for a decreasing sequence $\delta_n$ satisfying $(\delta_n)^\omega = o(n^{-1/2})$; (2) there exists $\Pi_n v^* \in \Gamma_n - \{\gamma_0\}$ such that $\delta_n \times \|\Pi_n v^* - v^*\| = o(n^{-1/2})$.

ASSUMPTION 4.6. For all $\tilde{\gamma} \in \mathcal{N}_{0n}$ with $\|\tilde{\gamma} - \gamma_0\| = O(\delta_n)$ and all $v = (v_\alpha, v_g)' \in \overline{\mathbf{V}}$ with $\|v\| = O(\delta_n)$ we have

$$E_0 \left( \frac{\partial^2 \ell(\tilde{\gamma}, Z_t)}{\partial \gamma \, \partial \gamma'}[v, v] - \frac{\partial^2 \ell(\gamma_0, Z_t)}{\partial \gamma \, \partial \gamma'}[v, v] \right) = o(n^{-1}).$$

For parametric likelihood models, Assumption 4.6 is automatically satisfied as long as the second-order derivatives of the log-likelihood are continuous in a shrinking neighborhood of the true parameter value. For sieve MLEs, Assumption 4.6 is satisfied provided that the third-order directional derivatives $\frac{d^3 \ell(\gamma_0 + \eta[\gamma - \gamma_0], Z_t)}{d\eta^3}$ exist for $\eta \in [0, 1]$, $\gamma \in \mathcal{N}_{0n}$ with $\|\gamma - \gamma_0\| = O(\delta_n)$, and the sieve MLE convergence rate $\delta_n$ is not too slow. For example, under Assumption 3.1(2) with polynomial, Fourier series, spline or wavelet sieves, we have a sieve MLE convergence rate of $\delta_n = n^{-r/(2r+1)}$ [see, e.g., Shen (1997) for i.i.d. data, and Chen and Shen (1998) for $\beta$-mixing time series data], and hence Assumption 4.6 is satisfied if $r > 1$.

ASSUMPTION 4.7. $\{\frac{\partial \ell(\gamma, Z_t)}{\partial \gamma'}[\Pi_n v^*] : \gamma \in \mathcal{N}_0, \|\gamma - \gamma_0\| = O(\delta_n)\}$ is a Donsker class.

Under Assumption 3.4(1), Assumption 4.7 is satisfied by applying the results of Doukhan, Massart and Rio (1995) on Donsker theorems for strictly stationary $\beta$-mixing processes.

THEOREM 4.1 (Normality). *Suppose that Assumptions M, 3.1–3.4 and 4.1–4.7 hold. Then* $\sqrt{n}(\rho(\widehat{\gamma}_n) - \rho(\gamma_0)) \Rightarrow N(0, \|\frac{\partial \rho(\gamma_0)}{\partial \gamma'}\|^2)$.

4.2. *Semiparametric efficiency of $\rho(\widehat{\gamma}_n)$.* We follow the approach of Wong (1992) to establish semiparametric efficiency. Related work can be found in Shen (1997), Bickel et al. (1993) and Bickel and Kwon (2001) and the references therein. Recall that a probability family $\{P_\gamma : \gamma \in \Gamma\}$ for the sample $\{Y_t\}_{t=1}^n$ is *locally asymptotically normal* (LAN) at $\gamma_0$, if (1) for any $v$ in the linear span of $\Gamma - \{\gamma_0\}$, $\gamma_0 + \eta n^{-1/2} v \in \Gamma$ for all small $\eta \geq 0$, and (2)

$$\frac{dP_{\gamma_0+n^{-1/2}v}}{dP_{\gamma_0}}(Y_1,\ldots,Y_n) = \exp\left\{n\left[L_n\left(\gamma_0 + \frac{1}{\sqrt{n}}v\right) - L_n(\gamma_0)\right]\right\}$$
$$= \exp\left\{\Sigma_n(v) - \frac{1}{2}\|v\|^2 + R_n(\gamma_0, v)\right\},$$

where $\Sigma_n(v)$ is linear in $v$, $\Sigma_n(v) \xrightarrow{d} \mathcal{N}(0, \|v\|^2)$ and $\operatorname{plim}_{n\to\infty} R_n(\gamma_0, v) = 0$ (both limits are under the true probability measure $P_{\gamma_0}$). To avoid the "super-efficiency" phenomenon, certain regularity conditions on the estimates are required. In estimating a smooth functional in the infinite-dimensional case, Wong [(1992), page 58] defines the class of *pathwise regular* estimates. An estimate $T_n(Y_1, \ldots, Y_n)$ of $\rho(\gamma_0)$ is *pathwise regular* if for any real number $\eta > 0$ and any $v$ in the linear span of $\Gamma - \{\gamma_0\}$, we have

$$\limsup_{n\to\infty} P_{\gamma_{n,\eta}}(T_n < \rho(\gamma_{n,\eta})) \leq \liminf_{n\to\infty} P_{\gamma_{n,-\eta}}(T_n < \rho(\gamma_{n,-\eta})),$$

where $\gamma_{n,\eta} = \gamma_0 + \eta n^{-1/2} v$ [see Wong (1992) and Shen (1997) for details].

THEOREM 4.2 (Efficiency). *Under conditions in Theorem 4.1, we have LAN, and the plug-in sieve MLE $\rho(\widehat{\gamma}_n)$ which achieves the efficiency lower bound for pathwise regular estimates.*

4.3. *$\sqrt{n}$ normality and efficiency of sieve MLE of copula parameter.* We take $\rho(\gamma) = \lambda'\alpha$ for any arbitrarily fixed $\lambda \in \mathcal{R}^d$ with $0 < |\lambda| < \infty$. It satisfies Assumption 4.4(2) with $\frac{\partial \rho(\gamma_0)}{\partial \gamma'}[v] = \lambda' v_\alpha$ and $\omega = \infty$. Assumption 4.4(1) is equivalent to finding a Riesz representer $v^* \in \overline{\mathbf{V}}$ satisfying (4.2) and (4.3),

$$(4.2) \qquad \lambda'(\alpha - \alpha_0) = \langle \gamma - \gamma_0, v^* \rangle \qquad \text{for any } \gamma - \gamma^* \in \overline{\mathbf{V}}$$

and

$$(4.3) \qquad \left\|\frac{\partial \rho(\gamma_0)}{\partial \gamma'}\right\|^2 = \|v^*\|^2 = \langle v^*, v^* \rangle = \sup_{v \neq 0, v \in \overline{\mathbf{V}}} \frac{|\lambda' v_\alpha|^2}{\|v\|^2} < \infty.$$

Let us change the variables before making statements on (4.3). Denote

$$\mathcal{L}_2^0([0,1]) \equiv \left\{e : [0,1] \to \mathcal{R} : \int_0^1 e(v)\,dv = 0, \int_0^1 [e(v)]^2\,dv < \infty\right\}.$$

By changing variables, for any $v_g \in \overline{\mathbf{V}}_g$, there is a unique function $b_g \in \mathcal{L}_2^0([0,1])$ with $b_g(u) = v_g(G_0^{-1}(u))/g_0(G_0^{-1}(u))$, and vice versa. So we can express $\frac{\partial \ell(\gamma_0, Z_t)}{\partial \gamma'}[v]$ as

$$\frac{\partial \ell(\gamma_0, Z_t)}{\partial \gamma'}[v] = \frac{\partial \ell(\gamma_0, U_t, U_{t-1})}{\partial \gamma'}[(v_\alpha', b_g)']$$

$$= \frac{\partial \log c(U_{t-1}, U_t; \alpha_0)}{\partial \alpha'}[v_\alpha] + b_g(U_t)$$

$$+ \sum_{j=1}^{2} \frac{\partial \log c(U_{t-1}, U_t; \alpha_0)}{\partial u_j} \int_0^{U_{t-2+j}} b_g(u)\, du$$

and

$$\|v\|^2 = E_0\left[\left(\frac{\partial \ell(\gamma_0, U_t, U_{t-1})}{\partial \gamma'}[(v_\alpha', b_g)']\right)^2\right]$$

$$= E_0\left[\left(\frac{\partial \log c(U_{t-1}, U_t; \alpha_0)}{\partial \alpha'}[v_\alpha] + b_g(U_t)\right.\right.$$

$$\left.\left. + \sum_{j=1}^{2} \frac{\partial \log c(U_{t-1}, U_t; \alpha_0)}{\partial u_j} \int_0^{U_{t-2+j}} b_g(u)\, du\right)^2\right].$$

Define

$$\overline{\mathbf{B}} = \left\{ b = (v_\alpha', b_g)' \in (\mathcal{A} - \alpha_0) \times \mathcal{L}_2^0([0,1]) : \|b\|^2 \right.$$

$$\left. \equiv E_0\left[\left(\frac{\partial \ell(\gamma_0, U_t, U_{t-1})}{\partial \gamma'}[b]\right)^2\right] < \infty \right\}.$$

Then there is a one-to-one onto mapping between the two Hilbert spaces $(\overline{\mathbf{B}}, \|\cdot\|)$ and $(\overline{\mathbf{V}}, \|\cdot\|)$. So the Riesz representer $v^* = (v_\alpha^{*\prime}, v_g^*)' \in \overline{\mathbf{V}}$ is uniquely determined by $b^* = (v_\alpha^{*\prime}, b_g^*)' \in \overline{\mathbf{B}}$ (and vice versa) via the relation $v_g^*(y) = b_g^*(G_0(y))g_0(y)$ for all $y \in \mathcal{Y}$. Notice that

$$\sup_{v \neq 0, v \in \overline{\mathbf{V}}} \frac{|\lambda' v_\alpha|^2}{\|v\|^2}$$

$$= \sup_{b \neq 0, b \in \overline{\mathbf{B}}} |\lambda' v_\alpha|^2 \bigg/ \left( E_0\left[\left(\frac{\partial \log c(U_{t-1}, U_t; \alpha_0)}{\partial \alpha'}[v_\alpha] + b_g(U_t)\right.\right.\right.$$

$$\left.\left. + \sum_{j=1}^{2} \frac{\partial \log c(U_{t-1}, U_t; \alpha_0)}{\partial u_j} \right.$$

$$\times \int_0^{U_{t-2+j}} b_g(u)\, du \bigg)^2 \bigg] \bigg)$$

$$= \lambda' \mathcal{I}_*(\alpha_0)^{-1} \lambda = \lambda' (E_0[\mathcal{S}_{\alpha_0} \mathcal{S}'_{\alpha_0}])^{-1} \lambda,$$

where $\mathcal{S}_{\alpha_0}$ is the efficient score function for $\alpha_0$,

$$
(4.4) \quad \begin{aligned}
\mathcal{S}'_{\alpha_0} &= \frac{\partial \log c(\alpha_0, U_t, U_{t-1})}{\partial \alpha'} - \mathbf{e}^*(U_t) \\
&\quad - \sum_{j=1}^{2} \frac{\partial \log c(\alpha_0, U_t, U_{t-1})}{\partial u_j} \int_0^{U_{t-2+j}} \mathbf{e}^*(u)\, du
\end{aligned}
$$

and $\mathbf{e}^* = (e_1^*, \ldots, e_d^*) \in (\mathcal{L}_2^0([0,1]))^d$ solves the following infinite-dimensional optimization problems for $k = 1, \ldots, d$,

$$\inf_{e_k \in \mathcal{L}_2^0([0,1])} E_0 \bigg\{ \bigg( \frac{\partial \log c(U_{t-1}, U_t; \alpha_0)}{\partial \alpha_k} - e_k(U_t) \\ - \sum_{j=1}^2 \frac{\partial \log c(U_{t-1}, U_t; \alpha_0)}{\partial u_j} \int_0^{U_{t-2+j}} e_k(u)\, du \bigg)^2 \bigg\}.$$

Therefore, $b^* = (v_\alpha^{*\prime}, b_g^*)'$ with $v_\alpha^* = \mathcal{I}_*(\alpha_0)^{-1} \lambda$ and $b_g^*(u) = -e^*(u) \times v_\alpha^*$, and $v^* = [I_d, -e^*(G_0(\cdot))g_0(\cdot)] \times \mathcal{I}_*(\alpha_0)^{-1} \lambda$. Hence (4.3) is satisfied if and only if $\mathcal{I}_*(\alpha_0) = E_0[\mathcal{S}_{\alpha_0} \mathcal{S}'_{\alpha_0}]$ is *nonsingular*, which in turn is satisfied under the following assumption:

ASSUMPTION 4.4'. (1) $\int \frac{\partial c(u;\alpha_0)}{\partial u_j}\, du_{-j} = \frac{\partial}{\partial u_j} \int c(u;\alpha_0)\, du_{-j} = 0$ for $(j, -j) = (1,2)$ with $j \neq -j$; (2) $\Sigma_{\text{ideal}} \equiv E_0(\frac{\partial \log c(U_{t-1}, U_t; \alpha_0)}{\partial \alpha} \{\frac{\partial \log c(U_{t-1}, U_t; \alpha_0)}{\partial \alpha}\}')$ is finite and positive definite; (3) $\int \frac{\partial^2 c(u;\alpha_0)}{\partial u_j \partial \alpha}\, du_{-j} = \frac{\partial^2}{\partial u_j \partial \alpha} \int c(u;\alpha_0)\, du_{-j} = 0$ for $(j, -j) = (1,2)$ with $j \neq -j$; (4) there exists a constant $K$ such that $\max_{j=1,2} \sup_{0 < u_j < 1} E[(u_j(1-u_j)\frac{\partial \log c(U_1, U_2; \alpha_0)}{\partial u_j})^2 | U_j = u_j] \leq K$.

Assumption 4.4' is a sufficient condition to ensure that the copula parameter can be estimated at a root-$n$ parametric rate. It is imposed in Bickel et al. (1993) and Chen, Fan and Tsyrennikov (2006) for semiparametric bivariate copula models. Bickel et al. (1993) has shown that many popular copula functions such as Clayton, Gaussian, Gumbel, Frank and others all satisfy this assumption. We can now apply Theorems 4.1 and 4.2 to obtain the following result:

PROPOSITION 4.1. *Suppose that Assumptions M, 3.1–3.4, 4.1–4.3, 4.4 and 4.5–4.7 hold. Then $\sqrt{n}(\widehat{\alpha}_n - \alpha_0) \Rightarrow N(0, \mathcal{I}_*(\alpha_0)^{-1})$, and $\widehat{\alpha}_n$ is semiparametrically efficient.*

In general, there is no closed-form solution of $\mathcal{I}_*(\alpha_0)$. Nevertheless, it can be consistently estimated by a sieve least squares method using its characterization in (4.4). Let $\widehat{U}_t = \widehat{G}_n(Y_t)$ for $t = 1, \ldots, n$. Let $\mathbf{B}_n$ be some sieve space such as

$$(4.5) \qquad \mathbf{B}_n = \left\{ e(u) = \sum_{k=1}^{K_{n\alpha}} a_k \sqrt{2} \cos(k\pi u), u \in [0,1], \sum_{k=1}^{K_{n\alpha}} a_k^2 < \infty \right\},$$

where $K_{n\alpha} \to \infty, (K_{n\alpha})^d/n \to 0$. For $k = 1, \ldots, d$, we compute $\widehat{e}_k$ as the solution to

$$\min_{e_k \in \mathbf{B}_n} \frac{1}{n-1} \sum_{t=2}^{n} \left( \frac{\partial \log c(\widehat{U}_{t-1}, \widehat{U}_t; \widehat{\alpha})}{\partial \alpha_k} - e_k(\widehat{U}_t) \right.$$

$$\left. - \sum_{j=1}^{2} \frac{\partial \log c(\widehat{U}_{t-1}, \widehat{U}_t; \widehat{\alpha})}{\partial u_j} \int_0^{\widehat{U}_{t-2+j}} e_k(u)\, du \right)^2.$$

Denote $\widehat{\mathbf{e}} = (\widehat{e}_1, \ldots, \widehat{e}_d)$ and

$$\widehat{\mathcal{I}}_* = \frac{1}{n-1} \sum_{t=2}^{n} \left\{ \left( \frac{\partial \log c(\widehat{U}_{t-1}, \widehat{U}_t; \widehat{\alpha})}{\partial \alpha'} - \widehat{\mathbf{e}}(\widehat{U}_t) \right.\right.$$

$$\left. - \sum_{j=1}^{2} \frac{\partial \log c(\widehat{U}_{t-1}, \widehat{U}_t; \widehat{\alpha})}{\partial u_j} \int_0^{\widehat{U}_{t-2+j}} \widehat{\mathbf{e}}(u)\, du \right)'$$

$$\times \left( \frac{\partial \log c(\widehat{U}_{t-1}, \widehat{U}_t; \widehat{\alpha})}{\partial \alpha'} - \widehat{\mathbf{e}}(\widehat{U}_t) \right.$$

$$\left.\left. - \sum_{j=1}^{2} \frac{\partial \log c(\widehat{U}_{t-1}, \widehat{U}_t; \widehat{\alpha})}{\partial u_j} \int_0^{\widehat{U}_{t-2+j}} \widehat{\mathbf{e}}(u)\, du \right) \right\}.$$

Following the proof of Theorem 5.1 in Ai and Chen (2003) we immediately obtain the following:

PROPOSITION 4.2. *Under all of the assumptions of Proposition 4.1, $\widehat{\mathcal{I}}_* = \mathcal{I}_*(\alpha_0) + o_p(1)$.*

4.4. *Sieve MLE of the marginal distribution.* Let us consider the estimation of $\rho(\gamma_0) = G_0(y)$ for some fixed $y \in \mathcal{Y}$ by the plug-in sieve MLE, $\rho(\widehat{\gamma}_n) = \widehat{G}_n(y) = \int 1(x \leq y) \widehat{g}_n(x)\, dx$, where $\widehat{g}_n$ is the sieve MLE for $g_0$.

Clearly, $\frac{\partial \rho(\gamma_0)}{\partial \gamma'}[v] = \int_{\mathcal{Y}} 1(x \leq y) v_g(x)\, dx$ for any $v = (v'_\alpha, v_g)' \in \overline{\mathbf{V}}$. It is easy to see that $\omega = \infty$ in Assumption 4.4, and

$$\left\| \frac{\partial \rho(\gamma_0)}{\partial \gamma'} \right\|^2 = \sup_{v \in \overline{\mathbf{V}}: \|v\|>0} \frac{|\int_{\mathcal{Y}} 1(x \leq y) v_g(x)\, dx|^2}{\|v\|^2} < \infty.$$

Hence the representer $v^* \in \overline{\mathbf{V}}$ should satisfy (4.6) and (4.7),

$$(4.6) \qquad \langle v^*, v \rangle = \frac{\partial \rho(\gamma_0)}{\partial \gamma'}[v] = E_0(1(Y_t \leq y)\frac{v_g(Y_t)}{g_0(Y_t)}) \qquad \text{for all } v \in \overline{\mathbf{V}},$$

$$(4.7) \qquad \left\|\frac{\partial \rho(\gamma_0)}{\partial \gamma'}\right\|^2 = \|v^*\|^2 = \|b^*\|^2 = \sup_{b \in \overline{\mathbf{B}}:\|b\|>0} \frac{|E_0[1(U_t \leq G_0(y))b_g(U_t)]|^2}{\|b\|^2}.$$

PROPOSITION 4.3. *Let $v^* \in \overline{\mathbf{V}}$ solve (4.6) and (4.7). Suppose that Assumptions M, 3.1–3.4, 4.1–4.3 and 4.5–4.7 hold. Then for any fixed $y \in \mathcal{Y}$, $\sqrt{n}(\widehat{G}_n(y) - G_0(y)) \Rightarrow N(0, \|v^*\|^2)$. Moreover, $\widehat{G}_n$ is semiparametrically efficient.*

Again, there are currently no closed-form expressions for the asymptotic variance $\|v^*\|^2$. Nevertheless, it can also be consistently estimated by the sieve method. Let

$$\widehat{\sigma}_G^2 \equiv \max_{v_\alpha \neq 0, b_g \in \mathbf{B}_n} \left|\frac{1}{n}\sum_{t=1}^n 1\{\widehat{U}_t \leq \widehat{G}_n(y)\}b_g(\widehat{U}_t)\right|^2$$

$$\bigg/ \left(\frac{1}{n-1}\sum_{t=2}^n \left[\frac{\partial \log c(\widehat{U}_{t-1}, \widehat{U}_t; \widehat{\alpha})}{\partial \alpha'}v_\alpha + b_g(\widehat{U}_t)\right.\right.$$

$$\left.\left. + \sum_{j=1}^2 \frac{\partial \log c(\widehat{U}_{t-1}, \widehat{U}_t; \widehat{\alpha})}{\partial u_j}\int_0^{\widehat{U}_{t-2+j}} b_g(u)\,du\right]^2\right),$$

where $\widehat{U}_t = \widehat{G}_n(Y_t)$, and $\mathbf{B}_n$ is given in (4.5).

PROPOSITION 4.4. *Under all the assumptions of Proposition 4.3, we have, for any fixed $y \in \mathcal{Y}$, $\widehat{\sigma}_G^2 = \|v^*\|^2 + o_p(1)$.*

4.5. *Plug-in estimates of conditional quantiles.* Under Assumption M, the $q$th conditional quantile of $Y_t$ given $Y_{t-1} = y$ is given by $Q_q^Y(y) = G_0^{-1}(C_{2|1}^{-1}[q|G_0(y); \alpha_0])$. Its plug-in sieve MLE estimate is given by

$$\widehat{Q}_q^Y(y) = \widehat{G}_n^{-1}(C_{2|1}^{-1}[q|\widehat{G}_n(y); \widehat{\alpha}_n]).$$

Let $\rho(\gamma_0) = Q_q^Y(y)$, then by some calculation, for any $v = (v_\alpha, v_g)' \in \overline{\mathbf{V}}$,

$$\frac{\partial \rho(\gamma_0)}{\partial \gamma'}[v] = \left(\frac{-C_{11}\int 1(x \leq y)v_g(x)\,dx - C_{1\alpha}v_\alpha}{c(U_{t-1}, C_1^{-1}(U_{t-1}, q; \alpha_0), \alpha_0)}\right.$$

$$\left. - \int 1(x \leq Q_q^Y(y))v_g(x)\,dx\right)\bigg/g_0(Q_q^Y(y))$$

where $C_{11} = \frac{\partial^2 C(U_{t-1}, C_1^{-1}(U_{t-1}, q; \alpha_0), \alpha_0)}{\partial u_1^2}$ and $C_{1\alpha} = \frac{\partial^2 C(U_{t-1}, C_1^{-1}(U_{t-1}, q; \alpha_0), \alpha_0)}{\partial u_1 \partial \alpha}$.

We can see $\omega = 2$ in Assumption 4.4, as long as $g_0(Q_q^Y(y)) \neq 0$ and $c(U_{t-1}, C_1^{-1}(U_{t-1}, q; \alpha_0), \alpha_0) \neq 0$, which are satisfied under Assumption M(2). Thus we have

$$\left\|\frac{\partial \rho(\gamma_0)}{\partial \gamma'}\right\|^2 = \sup_{v \in \overline{\mathbf{V}}: \|v\| > 0} \left| \{g_0(Q_q^Y(y))\}^{-1} \right.$$
$$\times \left[ \frac{-C_{11} \int 1(x \leq y) v_g(x)\, dx - C_{1\alpha} v_\alpha}{c(U_{t-1}, C_1^{-1}(U_{t-1}, q; \alpha_0), \alpha_0)} \right.$$
$$\left.\left. - \int 1(x \leq Q_q^Y(y)) v_g(x)\, dx \right]\right|^2 \Big/ \|v\|^2 < \infty.$$

Hence the Riesz representer $v^* \in \overline{\mathbf{V}}$ should satisfy: $\langle v^*, v \rangle = \frac{\partial \rho(\gamma_0)}{\partial \gamma'}[v]$ for all $v \in \overline{\mathbf{V}}$, and $\|v^*\|^2 = \|\frac{\partial \rho(\gamma_0)}{\partial \gamma'}\|^2$. Applying Theorems 4.1 and 4.2 we immediately obtain the following:

PROPOSITION 4.5. *Let $v^* \in \overline{\mathbf{V}}$ be the Riesz representer for $Q_q^Y(y)$. Suppose that Assumptions M, 3.1–3.4, 4.1–4.3 and 4.5–4.7 hold. Then for a fixed $y \in \mathcal{Y}$, $\sqrt{n}(\hat{Q}_q^Y(y) - Q_q^Y(y)) \Rightarrow N(0, \|v^*\|^2)$. Moreover, $\hat{Q}_q^Y(y)$ is semiparametrically efficient.*

**5. Sieve likelihood ratio inference for smooth functionals.** In this section, we are interested in the sieve likelihood ratio inference for smooth functional $\rho(\gamma) = (\rho_1(\gamma), \ldots, \rho_k(\gamma))' : \Gamma \to \mathcal{R}^k$,

$$H_0 : \rho(\gamma_0) = 0,$$

where $\rho$ is a vector of known functionals. [For instance, $\rho(\gamma) = \alpha - \alpha_0 \in \mathcal{R}^d$ or $\rho(\gamma) = G(y) - G_0(y) \in \mathcal{R}$ for fixed $y$.] Without loss of generality, we assume that $\frac{\partial \rho_1(\gamma_0)}{\partial \gamma'}, \ldots, \frac{\partial \rho_k(\gamma_0)}{\partial \gamma'}$ are linearly independent. Otherwise a linear transformation can be conducted for the hypothesis.

Suppose that $\rho_i$ satisfies Assumption 4.4 for $i = 1, \ldots, k$. Then by the Riesz representation theorem, there exists a $v_i^* \in \overline{\mathbf{V}}$ such that

$$\frac{\partial \rho_i(\gamma_0)}{\partial \gamma'}[v] \equiv \langle v_i^*, v \rangle \qquad \text{for all } v \in \overline{\mathbf{V}}.$$

Denote $v^* = (v_1^*, \ldots, v_k^*)'$. By the Gram–Schmidt orthogonalization, without loss of generality, we assume $\langle v_i^*, v_j^* \rangle = 0$ for any $i \neq j$.

Shen and Shi (2005) provide a theory on the sieve likelihood ratio inference for i.i.d. data. We now extend their result to strictly stationary Markov time series data. Denote

$$\hat{\gamma}_n = \arg \max_{\alpha \in \mathcal{A}, g \in \mathcal{G}_n} L_n(\alpha, g); \qquad \overline{\gamma}_n = \arg \max_{\alpha \in \mathcal{A}, g \in \mathcal{G}_n, \rho(\gamma) = 0} L_n(\alpha, g).$$

THEOREM 5.1. *Suppose that Assumptions M, 3.1–3.4, 4.1–4.3 and 4.5–4.7 hold, also that Assumption 4.4 holds with $\rho_i$, $i = 1, \ldots, k$, and Assumption 4.5(2) holds with $v_i^*$, $i = 1, \ldots, k$. Then*

$$2n(L_n(\widehat{\gamma}_n) - L_n(\overline{\gamma}_n)) \to^d \mathcal{X}^2_{(k)},$$

*where $\mathcal{X}^2_{(k)}$ stands for the chi-square distribution with $k$ degrees of freedom, and $\frac{\partial \rho_1(\gamma_0)}{\partial \gamma'}, \ldots, \frac{\partial \rho_k(\gamma_0)}{\partial \gamma'}$ are assumed to be linearly independent.*

We can apply Theorem 5.1 to construct confidence regions of any smooth functionals. For example, we can compute confidence region for sieve MLE of the copula parameter $\alpha$. Define $\tilde{g}_n(\alpha) = \arg\max_{g \in \mathcal{G}_n} L_n(\alpha, g)$. By Theorem 5.1, $2n(L_n(\widehat{\alpha}_n, \tilde{g}_n(\widehat{\alpha}_n)) - L_n(\alpha_0, \tilde{g}_n(\alpha_0))) \to^d \mathcal{X}^2_{(d)}$ where $(\widehat{\alpha}_n, \tilde{g}_n(\widehat{\alpha}_n)) = \widehat{\gamma}_n$ is the original sieve MLE.[4]

**6. Monte Carlo comparison of several estimators.** In this section, we address the finite sample performance of sieve MLE by comparing it to several existing popular estimators: the two-step semiparametric estimator proposed in Chen and Fan (2006), the ideal (or infeasible) MLE, the correctly specified parametric MLE and the misspecified parametric MLE.

6.1. *Existing estimators.* For comparison, we briefly review several existing estimators that have been used in applied work.

6.1.1. *Two-step semiparametric estimator.* Chen and Fan (2006) propose the following two-step semiparametric procedure:

*Step* 1. Estimate the unknown true marginal distribution $G_0(y)$ by the empirical distribution function $\frac{n+1}{n} G_n(y)$ where $G_n(y) \equiv \frac{1}{n+1} \sum_{t=1}^n 1\{Y_t \leq y\}$.

*Step* 2. Estimate the copula dependence parameter $\alpha_0$ by

$$\hat{\alpha}_n^{2sp} \equiv \arg\max_{\alpha \in \mathcal{A}} \frac{1}{n} \sum_{t=2}^n \log c(G_n(Y_{t-1}), G_n(Y_t); \alpha).$$

Assuming that the process $\{Y_t\}_{t=1}^n$ is $\beta$-mixing with a certain decay rate, under Assumption M and some other mild regularity conditions, Chen and Fan (2006) show that

$$\sqrt{n}(\hat{\alpha}_n^{2sp} - \alpha_0) \to_d N(0, \sigma_{2sp}^2), \qquad \text{with } \sigma_{2sp}^2 \equiv B_0^{-1} \Sigma_{2sp} B_0^{-1},$$

---

[4]If we only care about estimation and inference of copula parameter $\alpha$, we could also extend the results of Murphy and van der Vaart (2000) on profile likelihood ratio to our copula-based semiparametric Markov models.

where $B_0 \equiv -E_0(\frac{\partial^2 \log c(U_{t-1}, U_t; \alpha_0)}{\partial \alpha \partial \alpha'}) = \Sigma_{\text{ideal}}$ (under Assumption 4.4'), and

$$\Sigma_{2sp} \equiv \lim_{n \to \infty} \text{Var}_0 \left\{ \frac{1}{\sqrt{n}} \sum_{t=2}^n \left[ \frac{\partial \log c(U_{t-1}, U_t; \alpha_0)}{\partial \alpha} + W_1(U_{t-1}) + W_2(U_t) \right] \right\}$$
$$< \infty,$$

$$W_1(U_{t-1}) \equiv \int_0^1 \int_0^1 [1\{U_{t-1} \leq v_1\} - v_1] \frac{\partial^2 \log c(v_1, v_2; \alpha_0)}{\partial \alpha \, \partial u_1} c(v_1, v_2; \alpha_0) \, dv_1 \, dv_2,$$

$$W_2(U_t) \equiv \int_0^1 \int_0^1 [1\{U_t \leq v_2\} - v_2] \frac{\partial^2 \log c(v_1, v_2; \alpha_0)}{\partial \alpha \, \partial u_2} c(v_1, v_2; \alpha_0) \, dv_1 \, dv_2.$$

EXAMPLE 6.1 (Two-step semiparametric estimator of Gaussian copula parameter). The bivariate Gaussian copula is $C(u_1, u_2; \alpha) = \Phi_\alpha(\Phi^{-1}(u_1), \Phi^{-1}(u_2))$ for $|\alpha| < 1$ where $\Phi_\alpha$ is the bivariate standard normal distribution with correlation $\alpha$, and $\Phi$ is the scalar standard normal distribution. Chen and Fan (2006) show that

$$\sqrt{n}(\hat{\alpha}_n^{2sp} - \alpha_0) \to_d N(0, 1 - \alpha_0^2).$$

Klaassen and Wellner (1997) establish that the semiparametric efficient variance bound for estimating a Gaussian copula parameter $\alpha$ is $1 - \alpha_0^2$; hence $\hat{\alpha}_n^{2sp}$ is semiparametrically efficient for a Gaussian copula. However, as pointed out by Genest and Werker (2002), the Gaussian copula and the independence copula are the only two copulas for which the two-step semiparametric estimator is efficient for $\alpha_0$. Moreover, the empirical CDF estimator is still inefficient for $G_0(\cdot)$, even in this Gaussian copula-based Markov model.

6.1.2. *Possibly misspecified parametric MLE.* Denote $G(y, \theta)[g(y, \theta)]$ as the marginal distribution (marginal density) whose functional form is known up to the unknown finite-dimensional parameter $\theta$. Then the observed joint parametric log-likelihood for $\{Y_t\}_{t=1}^n$ is

$$L_n(\alpha, \theta) = \frac{1}{n} \sum_{t=1}^n \log g(Y_t, \theta) + \frac{1}{n} \sum_{t=2}^n \log c(G(Y_{t-1}, \theta), G(Y_t, \theta); \alpha),$$

and the parametric MLE is $(\hat{\alpha}_n^p, \hat{\theta}_n^p) = \arg\max_{(\alpha, \theta) \in \mathcal{A} \times \Theta} L_n(\alpha, \theta)$ where $\mathcal{A} \times \Theta$ is the parameter space. Under Assumption M and some other mild regularity conditions, we have

$$\sqrt{n}((\hat{\alpha}_n^p, \hat{\theta}_n^p) - (\alpha^*, \theta^*)) \to_d N(0, B_{*p}^{-1} \Sigma_{*p} B_{*p}^{-1}),$$

where $B_{*p} \equiv -E_0(\frac{\partial^2 \ell(\alpha^*, \theta^*, Z_t)}{\partial(\alpha, \theta) \, \partial(\alpha, \theta)'})$ is nonsingular and $\Sigma_{*p} \equiv \lim_{n \to \infty} \text{Var}\{\frac{1}{\sqrt{n}} \times \sum_{t=2}^n \frac{\partial \ell(\alpha^*, \theta^*, Z_t)}{\partial(\alpha, \theta)}\}$.

6.1.3. *Efficiency of correctly specified parametric MLE.* Asymptotic properties for correctly specified MLEs for Markov processes have been discussed in Section 10.4 of Joe (1997) and Billingsley (1961b). Under Assumption M and the correct specification of marginal $G(Y_t, \theta^*) = G_0(Y_t)$, we have $\alpha^* = \alpha_0$, $B_{*p} = \Sigma_{*p} = \Sigma_{0p} \equiv E_0(\frac{\partial \ell(\alpha_0, \theta^*, Z_t)}{\partial(\alpha, \theta)}\{\frac{\partial \ell(\alpha_0, \theta^*, Z_t)}{\partial(\alpha, \theta)}\}')$, and $(\hat{\alpha}_n^p, \hat{\theta}_n^p)$ is $\sqrt{n}$-efficient for $(\alpha_0, \theta^*)$ with asymptotic variance $\Sigma_{0p}^{-1}$. Moreover, $\sqrt{n}(\hat{\alpha}_n^p - \alpha_0) \to_d N(0, \mathcal{I}_{*p}(\alpha_0)^{-1})$ with

$$\mathcal{I}_{*p}(\alpha_0) \equiv \min_{\mathbf{b}} E_0\bigg(\bigg(\frac{\partial \log c(U_{t-1}, U_t; \alpha_0)}{\partial \alpha} - \frac{\partial \ell(\alpha_0, \theta^*, Z_t)}{\partial \theta}\mathbf{b}\bigg) \\ \times \bigg(\frac{\partial \log c(U_{t-1}, U_t; \alpha_0)}{\partial \alpha} - \frac{\partial \ell(\alpha_0, \theta^*, Z_t)}{\partial \theta}\mathbf{b}\bigg)'\bigg).$$

6.1.4. *Ideal (or infeasible) MLE.* We denote $\hat{\alpha}_n^{\text{Ideal}}$ as the ideal (or infeasible) MLE of the copula parameter $\alpha_0$ when the marginal $G_0(\cdot)$ is assumed to be completely known. Let $\hat{\alpha}_n^{\text{Ideal}} = \arg\max_{\alpha \in \mathcal{A}} \frac{1}{n} \sum_{t=2}^n \log c(U_{t-1}, U_t; \alpha)$. Suppose that Assumption M holds with a completely known $G(\cdot, \theta) = G_0(\cdot)$. Then $B_0 \equiv -E_0(\frac{\partial^2 \log c(U_{t-1}, U_t; \alpha_0)}{\partial \alpha \partial \alpha'}) = \Sigma_{\text{ideal}}$ is finite and nonsingular and $\hat{\alpha}_n^{\text{Ideal}}$ is efficient, thus

$$\sqrt{n}(\hat{\alpha}_n^{\text{Ideal}} - \alpha_0) \to_d N(0, \Sigma_{\text{ideal}}^{-1}).$$

REMARK 6.1. Since $\mathcal{I}_*(\alpha_0) \leq \mathcal{I}_{*p}(\alpha_0) \leq \Sigma_{\text{ideal}}$, we have $\mathcal{I}_*(\alpha_0)^{-1} \geq \mathcal{I}_{*p}(\alpha_0)^{-1} \geq \Sigma_{\text{ideal}}^{-1}$. Also, Proposition 4.1 immediately implies that $\sigma_{2sp}^2 \geq \mathcal{I}_*(\alpha_0)^{-1}$.

EXAMPLE 6.1' (The ideal MLE of a Gaussian copula parameter). For the Gaussian copula in Example 6.1, it is easy to verify that

$$\Sigma_{\text{ideal}} = B_0 = -E_0\bigg(\frac{\partial^2 \log c(U_{t-1}, U_t; \alpha_0)}{\partial \alpha \, \partial \alpha}\bigg) = \frac{1 + \alpha_0^2}{(1 - \alpha_0^2)^2} < \infty \quad \text{if } \alpha_0^2 \neq 1.$$

Consequently, $\sqrt{n}(\hat{\alpha}_n^{\text{Ideal}} - \alpha_0) \to_d N(0, \Sigma_{\text{ideal}}^{-1})$ with $\Sigma_{\text{ideal}}^{-1} = (1 - \alpha_0^2) \times \frac{1-\alpha_0^2}{1+\alpha_0^2}$. We note that the asymptotic variance $\text{Avar}(\hat{\alpha}_n^{\text{Ideal}}) = \Sigma_{\text{ideal}}^{-1} \leq 1 - \alpha_0^2 = \text{Avar}(\hat{\alpha}_n^{2sp})$, and $\text{Avar}(\hat{\alpha}_n^{\text{Ideal}}) = \text{Avar}(\hat{\alpha}_n^{2sp})$ if and only if $\alpha_0 = 0$ (i.e., independent copula). Also $\text{Avar}(\hat{\alpha}_n^{\text{Ideal}})$ is decreasing in $|\alpha_0|$.

EXAMPLE 2.1' (The ideal MLE of a Clayton copula parameter). For the Clayton copula in Example 2.1, after some tedious calculation, we have

$$\Sigma_{\text{ideal}} = B_0 = \frac{1}{\alpha(1+\alpha)} + \frac{1}{\alpha(1+\alpha)^2(1+2\alpha)} + \frac{(1+\alpha)(1+2\alpha)}{\alpha^5} \times \text{Int}(\alpha)$$

where $\text{Int}(\alpha) = \int_1^\infty \int_1^\infty \frac{xy(\log x - \log y)^2 - x(\log x)^2 - y(\log y)^2}{(x+y-1)^{4+1/\alpha}} \, dx \, dy$ which is a small number bounded in $[-1,1]$. Therefore, $\Sigma_{\text{ideal}} \in (0, \infty)$, provided that $\alpha_0 > 0$. Hence $\sqrt{n}(\hat{\alpha}_n^{\text{Ideal}} - \alpha_0) \to_d N(0, \Sigma_{\text{ideal}}^{-1})$ where the asymptotic variance $\Sigma_{\text{ideal}}^{-1}$ is increasing in $\alpha_0$ and is $O(\alpha_0^2)$.

EXAMPLE 2.1' (The ideal MLE of EFGM copula parameter). For the EFGM copula with $C(u_1, u_2; \alpha) = u_1 u_2 (1 + \alpha(1-u_1)(1-u_2))$, $\alpha \in [-1, 1]$, the copula density function is

$$c(u_1, u_2; \alpha) = \frac{\partial^2}{\partial u_1 \, \partial u_2} C(u_1, u_2; \alpha) = 1 + \alpha - 2\alpha(u_1 + u_2) + 4\alpha u_1 u_2.$$

Let $Li_2(z) = \sum_{k=1}^\infty z^k / k^2$, $|z| \leq 1$, be the polylogarithm function with order 2. Then

$$\Sigma_{\text{ideal}} = -E_0 \left( \frac{\partial^2 \log c(U_{t-1}, U_t; \alpha_0)}{\partial \alpha \, \partial \alpha} \right)$$
$$= \int_0^1 \int_0^1 \frac{(1 - 2u_1 - 2u_2 + 4u_1 u_2)^2}{1 + \alpha - 2\alpha(u_1 + u_2) + 4\alpha u_1 u_2} \, du_1 \, du_2$$
$$= \sum_{k=1}^\infty \frac{\alpha^{2k-2}}{(1+2k)^2} = \frac{Li_2(|\alpha|) - Li_2(\alpha^2)/4 - |\alpha|}{|\alpha|^3}.$$

6.2. *Simulations.* One can simulate a strictly stationary first-order Markov process $\{Y_t\}_{t=1}^n$ from a specified bivariate copula $C(u_1, u_2; \alpha_0)$ with given invariant CDF $G_0$ as follows:

*Step* 1. Generate an i.i.d. sequence of uniform random variables $\{V_t\}_{t=1}^n$.

*Step* 2. Set $U_1 = V_1$ and $U_t = C_{2|1}^{-1}[V_t | U_{t-1}, \alpha_0]$.

*Step* 3. Set $Y_t = G_0^{-1}(U_t)$ for $t = 1, \ldots, n$.

In our simulation study, we consider several first-order Markov models generated via different classes of copulas (Clayton, Gumbel, Frank, Gaussian and EFGM), with either Student $t_3$ or $t_5$ marginal distribution. Thus the true marginal distribution is $G_0 = t_\nu$ with density $g_0(y) = \frac{\Gamma(0.5(\nu+1))}{\sqrt{\nu\pi}\Gamma(\nu/2)} (1 + \frac{y^2}{\nu})^{-0.5(\nu+1)}$ with degrees of freedom $\nu = 3$ or 5. For each specified copula $C(u_1, u_2; \alpha_0)$, we generate a long time series, but we delete the first 2000 observations and keep the last 1000 observations as our simulated data sample data $\{Y_t\}$ (i.e., a simulated sample size $n = 1000$).

For all the copula-based Markov models and for each simulated sample, we compute five estimators of $\alpha_0$: sieve MLEs, ideal (or infeasible) MLEs, two-step estimators, correctly specified parametric MLEs (when the functional form of $g$ is correctly specified) and misspecified parametric MLEs (when the functional form of $g$ is misspecified). Sieve MLEs are computed by maximizing the joint log-likelihood $L_n(\alpha, g)$ in (3.1) using either a polynomial sieve or a polynomial spline sieve to approximate the log-marginal density ($\log g$). The selection of $K$, the number of the sieve terms, is based on the so-called small sample AIC of Burnham and Anderson (2002), $\widehat{K} = \arg\max_K \{L_n(\widehat{\gamma}_n(K)) - K/(n-K-1)\}$, where $\widehat{\gamma}_n(K)$ is the sieve MLE of $\gamma_0 = (\alpha_0, g_0)$ using $K$ as the sieve number of terms.

We compare the estimates of the copula parameter, and the estimates of $1/3$ and $2/3$ of the marginal quantiles in terms of Monte Carlo means, biases, variances, mean squared errors and confidence regions based on 1000 Monte Carlo simulation runs.

*Brief summary of MC results.* In the longer version posted on arXiv [Chen, Wu and Yi (2009)], we report all the simulation findings in detail. Here we only report a few Monte Carlo results for Clayton and Gumbel copula-based Markov models in Appendix B, and give a brief summary of the overall patterns. (1) Sieve MLEs of copula parameters always perform better than the two-step estimator in terms of bias and MSE, except for Gaussian copulas and EFGM copulas. For Gaussian copulas, we already explained (in Example 6.1) that both the sieve MLE and the two-step estimators are semiparametrically efficient for the copula parameter with unknown marginal distributions. For EFGM copulas, the distance between the EFGM copula function and the independent copula function is $\alpha u_1 u_2 (1-u_1)(1-u_2) \leq 0.0625\alpha$ for $\alpha \in [-1,1]$. Therefore, the EFGM copula is very close to the independent copula; hence the performance of the sieve MLE, the two-step, the correctly specified parametric MLE and the ideal MLE for copula parameters are all very close to one another; (2) For all the copula-based Markov models with some dependence in terms of Kendall's $\tau \neq 0$, including Gaussian and EFGM copula-based Markov models, sieve MLEs of marginal distributions always perform better than the empirical CDFs in terms of bias and MSE; (3) For Markov models generated via strong tail dependent copulas, both the two-step-based estimators of copula parameters and the empirical CDF estimator of the marginal distribution perform very poorly, both having big biases and big MSEs; (4) Sieve MLEs perform very well even for copulas with strong tail dependence and fat-tailed marginal density $t_3$; (5) Extreme conditional quantiles estimated via sieve MLEs are much more precise than those estimated via two-step estimators; (6) Misspecified parametric MLEs could lead to inconsistent estimation of copula parameters (in addition to inconsistent estimation of marginal density parameters). In summary we recommend sieve MLEs to estimate copula-based Markov models and their implied conditional quantiles (VaRs).

**7. Conclusions.** In this paper, we first show that several widely used tail dependent copula-generated Markov models are in fact geometrically ergodic (hence geometrically $\beta$-mixing), although their time series plots may look highly persistent and "long memory alike." We then propose sieve MLEs for the class of first-order strictly stationary copula-based semiparametric Markov models that are characterized by the parametric copula parameter $\alpha_0$ and the unknown invariant density $g_0(\cdot)$. We show that the sieve MLEs of any smooth functional of $(\alpha_0, g_0)$ is root-$n$ consistent, asymptotically normal and efficient, and that their sieve likelihood ratio statistics are asymptotically chi-square distributed. We propose either consistent plug-in estimation of the asymptotic variance or inverting the sieve likelihood ratio statistics to construct confidence regions for the sieve MLEs. Monte Carlo studies indicate that, even for semiparametric Markov models generated via tail dependent copulas with fat-tailed marginal distributions, the sieve MLEs of the copula parameter, the marginal CDFs and the conditional quantiles all perform very well in finite samples.

In this paper, we assume that the parametric copula function is correctly specified. We could test this assumption by performing a sieve likelihood ratio test [see e.g., Fan and Jiang (2007) for a review about generalized likelihood ratio tests]. Alternatively, we could also consider a joint sieve ML estimation of nonparametric copulas and nonparametric marginals. Recently, Chen, Peng and Zhao (2009) provided an empirical likelihood estimation of nonparametric copulas using a bivariate random sample; their method could be extended to our time series setting.

## APPENDIX A: MATHEMATICAL PROOFS

We first recall some equivalent definitions of $\beta$-mixing and ergodicity for strictly stationary Markov processes. Then we present the drift criterion for geometric ergodicity of Markov chains.

DEFINITION A.1. (1) [Davydov (1973)] For a strictly stationary Markov process $\{Y_t\}_{t=1}^\infty$, the $\beta$-mixing coefficients are given by

$$\beta_t = \int \sup_{0 \leq \phi \leq 1} |E[\phi(Y_{t+1})|Y_1 = y] - E[\phi(Y_{t+1})]| \, dG_0(y).$$

The process $\{Y_t\}$ is $\beta$-mixing if $\lim_{t \to \infty} \beta_t = 0$; and it is geometric $\beta$-mixing if $\beta_t \leq \gamma \exp(-\delta t)$ for some $\delta, \gamma > 0$.

(2) [Chan and Tong (2001)] A strictly stationary Markov process $\{Y_t\}$ is (Harris) ergodic if

$$\lim_{t \to \infty} \sup_{0 \leq \phi \leq 1} |E[\phi(Y_{t+1})|Y_1 = y] - E[\phi(Y_{t+1})]| = 0 \qquad \text{for almost all } y;$$

and it is geometrically ergodic if there existsa measurable function $W$ with $\int W(y) \, dG_0(y) < \infty$ and a constant $\kappa \in [0, 1)$ such that for all $t \geq 1$,

(A.1) $$\sup_{0 \leq \phi \leq 1} |E[\phi(Y_{t+1})|Y_1 = y] - E[\phi(Y_{t+1})]| \leq \kappa^t W(y).$$

DEFINITION A.2. Let $\{Y_t\}$ be an irreducible Markov chain with a transition measure $P^n(y; A) = P(Y_{t+n} \in A | Y_t = y)$, $n \geq 1$. A nonnull set $S$ is called small if there exists a positive integer $n$, a constant $b > 0$ and a probability measure $\nu(\cdot)$ such that $P^n(y; A) \geq b\nu(A)$ for all $y \in S$ and all measurable set $A$s.

THEOREM A.1 [Theorem B.1.4 in Chan and Tong (2001)]. *Let $\{Y_t\}$ be an irreducible and aperiodic Markov chain. Suppose there exists a small set $S$, a nonnegative measurable function $L$ which is bounded away from 0 and $\infty$ on $S$, and constants $r > 1$, $\gamma > 0, K > 0$ such that*

$$(A.2) \qquad rE[L(Y_{t+1})|Y_t = y] \leq L(y) - \gamma \qquad \text{for all } y \notin S,$$

*and, let $S'$ be the complement of $S$,*

$$(A.3) \qquad \int_{S'} L(w) P(y, dw) < K \qquad \text{for all } y \in S.$$

*Then $\{Y_t\}$ is geometrically ergodic and (A.1) holds. Here $L$ is called the Lyapunov function.*

PROOF OF THEOREM 2.1. We establish the results by applying Theorem A.1 or applying Proposition 2.1(i) of Chen and Fan (2006).

(1) For a Clayton copula, let $\{Y_t\}_{t=1}^n$ be a stationary Markov process of order 1 generated from a bivariate Clayton copula and a marginal CDF $G_0(\cdot)$. Then the transformed process $\{U_t \equiv G_0(Y_t)\}_{t=1}^n$ has uniform marginals and a Clayton copula joint distribution of $(U_{t-1}, U_t)$. When $\alpha = 0$, the Clayton copula becomes the independence copula; hence the process $\{U_t \equiv G_0(Y_t)\}_{t=1}^n$ is i.i.d. and trivially geometrically ergodic.

Let $\alpha > 0$. Recall that $C_{2|1}[w|u;\alpha] = \frac{\partial}{\partial u} C(u, w; \alpha) = (u^{-\alpha} + w^{-\alpha} - 1)^{-1-1/\alpha} u^{-1-\alpha}$ and that $C_{2|1}^{-1}[q|u;\alpha_0] = [(q^{-\alpha/(1+\alpha)} - 1)u^{-\alpha} + 1]^{-1/\alpha}$ is the $q$th conditional quantile of $U_t$ given $U_{t-1} = u$. Denote $X_t \equiv U_t^{-\alpha}$. Let $\{V_t\}_{t=1}^n$ be a sequence of i.i.d. uniform$(0, 1)$ random variables such that $V_t$ is independent of $U_{t-1}$. Let $q = V_t$ in the above conditional quantile expression of $U_t$ given $U_{t-1}$, then we obtain the following nonlinear AR(1) model from the Clayton copula:

$$X_t = (V_t^{-\alpha/(1+\alpha)} - 1)X_{t-1} + 1 \qquad \text{with } X_t^{-1/\alpha} \equiv U_t \sim \text{uniform}(0, 1).$$

Note that the state space of $\{X_t\}$ is $(1, \infty)$. Since $E_0[(V_t^{-\alpha/(1+\alpha)} - 1)^{1/\alpha}] = 1$, we can let $p \in (0, 1/\alpha)$, and $L(x) = x^p > 1$ be the Lyapunov function. Then by Hölder's inequality, $\rho \equiv E_0[L(V_t^{-\alpha/(1+\alpha)} - 1)] < 1$. Let $r = \rho^{-1/2} > 1$ and

$$x_0 = \max\{x \geq 1 : rE_0[|x(V_t^{-\alpha/(1+\alpha)} - 1) + 1|^p] \geq x^p - 1\}.$$

Such $x_0$ always exists since

$$\lim_{x\to\infty} \frac{rE_0[|x(V_t^{-\alpha/(1+\alpha)} - 1) + 1|^p]}{x^p - 1} = r\rho = \rho^{1/2} < 1.$$

Let set $S = [1, x_0]$. Clearly, $L$ is bounded away from 0 and $\infty$ on $S$. We now show that $S$ is a small set. Let $f(\cdot|x)$ be the conditional density function of $X_1$ given $X_0 = x$. Then

$$f(y|x) = \frac{(1+\alpha)x^{1+1/\alpha}}{\alpha(y-1+x)^{2+1/\alpha}} \geq \frac{1+\alpha}{\alpha(y-1+x_0)^{2+1/\alpha}}$$

if $x \in S$. Choose the probability measure $\nu$ on $(1, \infty)$ as $\nu(dy) = f(y|x_0)\,dy$. Then

$$\Pr(X_1 \in A | X_0 = x) \geq x_0^{-1-1/\alpha} \nu(A) \quad \text{for all } x \in S \text{ and } A \in \mathcal{B}.$$

Hence $S$ is a small set; see Definition A.2. Notice that, by the definition of $x_0$,

$$rE_0[L(X_1)|X_0 = x] \leq L(x) - 1 \quad \text{for all } x > x_0,$$
$$E_0[L(X_1)|X_0 = x] < \infty \quad \text{for all } x \in S = [1, x_0].$$

Thus all of the conditions in Theorem A.1 are satisfied; hence $\{X_t\}_{t=1}^n$ is geometrically ergodic, and geometric $\beta$-mixing.

(2) For the Gumbel copula, let $\{Y_t\}_{t=1}^n$ be a stationary Markov process of order 1 generated from a bivariate Gumbel copula and a marginal CDF $G_0(\cdot)$. Then the transformed process $\{U_t \equiv G_0(Y_t)\}_{t=1}^n$ has uniform marginals and $(U_{t-1}, U_t)$ has the Gumbel copula joint distribution (see Example 2.2). When $\alpha = 1$, the Gumbel copula becomes the independence copula; hence the process $\{U_t \equiv G_0(Y_t)\}_{t=1}^n$ is i.i.d. and trivially geometrically ergodic.

Let $\alpha > 1$. Let $X_t = (-\log U_t)^\alpha$. Then $U_t = F(X_t)$, with $F(x) = \exp\{-x^{1/\alpha}\}$. Let $f(x) = \alpha^{-1} x^{1/\alpha - 1} \exp\{-x^{1/\alpha}\}$. Then for $X_t$ we have

$$\Pr(X_{t+1} \geq x_2 | X_t = x_1) = \frac{f(x_1 + x_2)}{f(x_1)}, \quad x_1, x_2 > 0.$$

Hence

$$E_0(X_{t+1}|X_t = x_1) = \int_0^\infty \Pr(X_{t+1} \geq x_2|X_t = x_1)\,dx_2 = \int_0^\infty \frac{f(x_1 + x_2)}{f(x_1)}\,dx_2$$
$$= \frac{F(x_1)}{f(x_1)} = \alpha x_1^{1-(1/\alpha)}.$$

Note that as $x_1 \to 0$,

$$E_0(X_{t+1}^{-1/(2\alpha)}|X_t = x_1) = \int_0^\infty x_2^{-1/(2\alpha)} \frac{-f'(x_1 + x_2)}{f(x_1)} dx_2$$

$$= x_1^{1-1/(2\alpha)} \int_0^\infty u^{-1/(2\alpha)} \frac{-f'(x_1 + x_1 u)}{f(x_1)} du$$

$$\sim x_1^{-1/(2\alpha)}(1 - 1/\alpha) \int_0^1 t^{-1/(2\alpha)}(1-t)^{-1/(2\alpha)} dt$$

where the last relation is due to

$$\lim_{x_1 \to 0} \frac{-f'(x_1 + x_1 u)}{f(x_1)} \times x_1 = (1 - 1/\alpha)(1+u)^{1/\alpha - 2}.$$

Observe that, as $\alpha > 1$,

$$\kappa_\alpha \equiv (1 - 1/\alpha) \int_0^1 t^{-1/(2\alpha)}(1-t)^{-1/(2\alpha)} dt$$

$$= (1 - 1/\alpha) \times B(1 - 1/(2\alpha), 1 - 1/(2\alpha)) < 1,$$

where $B(\cdot, \cdot)$ is the beta function.

Let $L(x) = x^{-1/(2\alpha)} + x$ be the Lyapunov function. Let $r = \inf_{x>0} L(x)/2$. Then

$$\lim_{x \to \infty} \frac{E_0(L(X_{t+1})|X_t = x)}{L(x) - r} = 0$$

and

$$\lim_{x \to 0} \frac{E_0(L(X_{t+1})|X_t = x)}{L(x) - r} = \kappa_\alpha < 1.$$

Let $S = [1/\lambda, \lambda]$ with sufficiently large $\lambda > 0$. Then $S$ is a small set. So all conditions in Theorem A.1 are satisfied; hence $\{X_t\}_{t=1}^n$ is geometrically ergodic and geometrically $\beta$-mixing.

(3) For Student's $t$ copula, let $\{Y_t\}_{t=1}^n$ be a stationary Markov process of order 1 generated from a bivariate $t$-copula and a marginal CDF $G_0(\cdot)$. Then the transformed process $\{U_t \equiv G_0(Y_t)\}_{t=1}^n$ satisfies the following:

$$t_\nu^{-1}(U_t) = \rho t_\nu^{-1}(U_{t-1}) + e_t \sqrt{\frac{\nu + (t_\nu^{-1}(U_{t-1}))^2}{\nu + 1}(1 - \rho^2)},$$

where $e_t \sim t_{\nu+1}$, and is independent of $U^{t-1} \equiv (U_{t-1}, \ldots, U_1)$ [see, e.g., Chen, Koenker and Xiao (2008)]. Let $X_t \equiv t_\nu^{-1}(U_t)$. Then

$$X_t = \rho X_{t-1} + \sigma(X_{t-1})e_t, \qquad \sigma(X_{t-1}) = \sqrt{\frac{\nu + (X_{t-1})^2}{\nu + 1}(1 - \rho^2)},$$

where $e_t \sim t_{\nu+1}$, and is independent of $X^{t-1} \equiv (X_{t-1}, \ldots, X_1)$. Let $L(x) = |x| + 1 \geq 1$ be the Lyapunov function. Then $E_0\{L(X_t)\} = \sqrt{\nu}\frac{\Gamma(\frac{\nu-1}{2})}{\sqrt{\pi}\Gamma(\nu/2)} + 1 < \infty$ provided that $\nu > 1$. Then

$$\begin{aligned}E_0(L(X_t)|X_{t-1}=x) &= E_0(|\rho X_{t-1} + \sigma(X_{t-1})e_t| \mid X_{t-1}=x) + 1 \\ &= E_0(|\rho x + \sigma(x)e_t|) + 1 \\ &< \sqrt{E_0(|\rho x + \sigma(x)e_t|^2)} + 1 \\ &= \sqrt{(\rho^2 x^2 + \sigma^2(x)E_0[e_t^2])} + 1,\end{aligned}$$

where the strict inequality is due to $e_t \sim t_{\nu+1}$ and for fixed $x$,

$$0 < \mathrm{Var}(|\rho x + \sigma(x)e_t|^2) = E(|\rho x + \sigma(x)e_t|^2) - [E_0(|\rho x + \sigma(x)e_t|)]^2.$$

Since $\sigma^2(x) = (1-\rho^2)(\nu + x^2)/(\nu+1)$, we have

$$\begin{aligned}\lim_{|x|\to\infty} &\frac{E_0(L(X_t)|X_{t-1}=x)}{L(x)} \\ &= \lim_{|x|\to\infty} \frac{E_0(|\rho x + \sigma(x)e_t|) + 1}{|x|+1} \\ &< \lim_{|x|\to\infty} \frac{\sqrt{(\rho^2 x^2 + \sigma^2(x)E_0[e_t^2])} + 1}{|x|+1} \\ &= \sqrt{\rho^2 + \frac{1-\rho^2}{\nu+1}E_0[e_t^2]} \leq \sqrt{\rho^2 + \frac{1-\rho^2}{2+1}E[t_3^2]} = 1,\end{aligned}$$

where the last inequality is due to $E_0[e_t^2]/(\nu+1)$ decreasing in $\nu \in [2,\infty]$, and the last equality is due to $E[t_3^2] = 3$. Then we can choose a small set $S = [-x_0, x_0]$ with sufficiently large $x_0 > 0$. Clearly the density of $e_t$ is bounded from above and below on a compact set. Hence, all conditions in Theorem A.1 or in Proposition 2.1(i) of Chen and Fan (2006) are satisfied, and $\{X_t\}_{t=1}^n$ is geometrically ergodic (hence geometrically $\beta$-mixing). $\square$

PROOF OF PROPOSITION 3.1. Since most of the conditions of consistency in Theorem 3.1 of Chen (2007) are already assumed in our Assumptions M, 3.1 and 3.2, it suffices to verify Condition 3.5 (uniform convergence over sieves) of Chen (2007). Assumption M implies that $\{Y_t\}_{t=1}^n$ is stationary ergodic. This and Assumption 3.2 implies that Glivenko–Cantelli theorem for a stationary ergodic processes is applicable, and hence $\sup_{\gamma \in \Gamma_n} |L_n(\gamma) - E\{L_n(\gamma)\}| = o_p(1)$. The result now follows from Theorem 3.1 of Chen (2007). $\square$

PROOF OF LEMMA 4.1. For (1), recall that $Z_t = (Y_{t-1}, Y_t)$, under Assumptions M, 3.1(1), (2), 4.1 and 4.2, we have, for all $s < t$,

$$E_0\bigg(\bigg(\frac{\partial \ell(\gamma_0, Z_t)}{\partial \gamma'}[v]\bigg)\bigg(\frac{\partial \ell(\gamma_0, Z_s)}{\partial \gamma'}[\tilde{v}]\bigg)\bigg)$$

$$= E_0\bigg(E_0\bigg(\bigg(\frac{\partial \ell(\gamma_0, Z_t)}{\partial \gamma'}[v]\bigg)\bigg(\frac{\partial \ell(\gamma_0, Z_s)}{\partial \gamma'}[\tilde{v}]\bigg)\bigg|Y_1,\ldots,Y_{t-1}\bigg)\bigg)$$

$$= E_0\bigg(\bigg(\frac{\partial \ell(\gamma_0, Z_s)}{\partial \gamma'}[\tilde{v}]\bigg)E_0\bigg(\frac{\partial \ell(\gamma_0, Z_t)}{\partial \gamma'}[v]\bigg|Y_{t-1}\bigg)\bigg).$$

Recall that the true conditional density function is $p^0(Y_t|Y^{t-1}) = g_0(Y_t) \times c(G_0(Y_{t-1}), G_0(Y_t); \alpha_0) = h(Y_t|Y_{t-1}; \gamma_0)$. We have

$$E_0\bigg(\frac{\partial \ell(\gamma_0, Z_t)}{\partial \gamma'}[v]\bigg|Y_{t-1}\bigg)$$

$$= \int \frac{\frac{\partial h(y_t|Y_{t-1};\gamma_0)}{\partial \gamma'}}{h(y_t|Y_{t-1};\gamma_0)}[v] h(y_t|Y_{t-1};\gamma_0)\, dy_t$$

$$= \int \frac{\partial h(y_t|Y_{t-1};\gamma_0)}{\partial \gamma'}[v]\, dy_t$$

$$= \frac{d(\int h(y_t|Y_{t-1};\gamma_0 + \eta v)\, dy_t)}{d\eta}\bigg|_{\eta=0} = \frac{d(1)}{d\eta}\bigg|_{\eta=0} = 0,$$

where the order of differentiation and integration can be reversed due to Assumption 4.3.

For (2), the above equality also implies that $\{\frac{\partial \ell(\gamma_0, Z_t)}{\partial \gamma'}[v]\}_{t=1}^n$ is a martingale difference sequence with respect to the filtration $\mathcal{F}_{t-1} = \sigma(Y_1; \ldots; Y_{t-1})$.

For (3), Since $\int h(y|Y_{t-1};\gamma_0 + \eta v)\, dy \equiv 1$, by differentiating this equation with respect to $\eta$ twice and evaluating it at $\eta = 0$, we get $E_0((\frac{\partial \ell(\gamma_0, Z_t)}{\partial \gamma'}[v])^2|Y_{t-1}) = -E_0(\frac{\partial^2 \ell(\gamma_0, Z_t)}{\partial \gamma \partial \gamma'}[v, v]|Y_{t-1})$ where the interchange of differentiation and integration is guaranteed by Assumption 4.3. □

PROOF OF THEOREM 4.1. Let $\epsilon_n$ be any positive sequence satisfying $\epsilon_n = o(n^{-1/2})$. Denote $r[\gamma, \gamma_0, Z_t] \equiv \ell(\gamma, Z_t) - \ell(\gamma_0, Z_t) - \frac{\partial \ell(\gamma_0, Z_t)}{\partial \gamma'}[\gamma - \gamma_0]$ and $\mu_n(g(Z_t)) = n^{-1}\sum_{t=2}^n [g(Z_t) - E_0 g(Z_t)]$. In the proof we let $g(\cdot)$ be $\ell(\gamma, \cdot)$, $r[\gamma, \gamma_0, \cdot]$ or $\frac{\partial \ell(\gamma_0, \cdot)}{\partial \gamma'}[v^*]$. Then by the definition of the sieve MLE $\widehat{\gamma}_n$ (with abuse of notation, we denote it as $\hat{\gamma}$ in the following),

$$0 \leq \frac{1}{n}\sum_{t=2}^n [\ell(\hat{\gamma}, Z_t) - \ell(\hat{\gamma} \pm \epsilon_n \Pi_n v^*, Z_t)]$$

$$= \mu_n(\ell(\hat{\gamma}, Z_t) - \ell(\hat{\gamma} \pm \epsilon_n \Pi_n v^*, Z_t))$$

$$+ E_0(\ell(\hat{\gamma}, Z_t) - \ell(\hat{\gamma} \pm \epsilon_n \Pi_n v^*, Z_t)) + o_p(n^{-1})$$

$$= \mp \epsilon_n \frac{1}{n} \sum_{t=2}^{n} \frac{\partial \ell(\gamma_0, Z_t)}{\partial \gamma'} [\Pi_n v^*] + \mu_n(r[\hat{\gamma}, \gamma_0, Z_t] - r[\hat{\gamma} \pm \epsilon_n \Pi_n v^*, \gamma_0, Z_t])$$

$$+ E_0(r[\hat{\gamma}, \gamma_0, Z_t] - r[\hat{\gamma} \pm \epsilon_n \Pi_n v^*, \gamma_0, Z_t]) + o(n^{-1}).$$

CLAIM 1. $\frac{1}{n} \sum_{t=2}^{n} \frac{\partial \ell(\gamma_0, Z_t)}{\partial \gamma'} [\Pi_n v^* - v^*] = o_p(n^{-1/2})$. This claim is true due to Chebyshev's inequality, serially uncorrelated (Lemma 4.1) and identically distributed data, and $\|\Pi_n v^* - v^*\| = o(1)$.

CLAIM 2. $\mu_n(r[\hat{\gamma}, \gamma_0, Z_t] - r[\hat{\gamma} \pm \epsilon_n \Pi_n v^*, \gamma_0, Z_t]) = \epsilon_n \times o_p(n^{-1/2})$. This claim holds since

$$\mu_n(r[\hat{\gamma}, \gamma_0, Z_t] - r[\hat{\gamma} \pm \epsilon_n \Pi_n v^*, \gamma_0, Z_t])$$

$$= \mu_n \left( \ell(\hat{\gamma}, Z_t) - \ell(\hat{\gamma} \pm \epsilon_n \Pi_n v^*, Z_t) \pm \epsilon_n \frac{\partial \ell(\gamma_0, Z_t)}{\partial \gamma'} [\Pi_n v^*] \right)$$

$$= \mp \epsilon_n \mu_n \left( \frac{\partial \ell(\tilde{\gamma}, Z_t)}{\partial \gamma'} [\Pi_n v^*] - \frac{\partial \ell(\gamma_0, Z_t)}{\partial \gamma'} [\Pi_n v^*] \right) = \epsilon_n \times o_p(n^{-1/2}),$$

where $\tilde{\gamma} \in \Gamma_n$ lies between $\hat{\gamma}$ and $\hat{\gamma} \pm \epsilon_n \Pi_n v^*$, and the last equality is implied by Assumption 4.7.

CLAIM 3. $E_0(r[\hat{\gamma}, \gamma_0, Z_t] - r[\hat{\gamma} \pm \epsilon_n \Pi_n v^*, \gamma_0, Z_t]) = \pm \epsilon_n \langle \hat{\gamma} - \gamma_0, v^* \rangle + \epsilon_n o_p(n^{-1/2}) + o_p(n^{-1})$.

Note that

$$E_0(r[\gamma, \gamma_0, Z_t]) = E_0 \left( \ell(\gamma, Z_t) - \ell(\gamma_0, Z_t) - \frac{\partial \ell(\gamma_0, Z_t)}{\partial \gamma'} [\gamma - \gamma_0] \right)$$

$$= \frac{1}{2} E_0 \left( \frac{\partial^2 \ell(\tilde{\gamma}, Z_t)}{\partial \gamma \partial \gamma'} [\gamma - \gamma_0, \gamma - \gamma_0] \right.$$

$$\left. - \frac{\partial^2 \ell(\gamma_0, Z_t)}{\partial \gamma \partial \gamma'} [\gamma - \gamma_0, \gamma - \gamma_0] \right)$$

$$+ \frac{1}{2} E_0 \left( \frac{\partial^2 \ell(\gamma_0, Z_t)}{\partial \gamma \partial \gamma'} [\gamma - \gamma_0, \gamma - \gamma_0] \right) + \epsilon_n \times o_p(n^{-1/2})$$

$$= \frac{1}{2} E_0 \left( \frac{\partial^2 \ell(\gamma_0, Z_t)}{\partial \gamma \partial \gamma'} [\gamma - \gamma_0, \gamma - \gamma_0] \right)$$

$$+ \epsilon_n \times o_p(n^{-1/2}) + o_p(n^{-1}),$$

where $\tilde{\gamma} \in \Gamma_n$ is located between $\gamma$ and $\gamma_0$, and the last equality is due to Assumption 4.6. By Lemma 4.1(3), we have

$$\|\gamma - \gamma_0\|^2$$
$$\equiv E_0\left[\left(\frac{\partial \ell(\gamma_0, Z_t)}{\partial \gamma'}[\gamma - \gamma_0]\right)^2\right] = -E_0\left(\frac{\partial^2 \ell(\gamma_0, Z_t)}{\partial \gamma \partial \gamma'}[\gamma - \gamma_0, \gamma - \gamma_0]\right).$$

Therefore,

$$E_0(r[\hat{\gamma}, \gamma_0, Z_t] - r[\hat{\gamma} \pm \epsilon_n \Pi_n v^*, \gamma_0, Z_t])$$
$$= -\frac{\|\hat{\gamma} - \gamma_0\|^2 - \|\hat{\gamma} \pm \epsilon_n \Pi_n v^* - \gamma_0\|^2}{2} + o_p(\epsilon_n n^{-1/2}) + o_p(n^{-1})$$
$$= \pm \epsilon_n \langle \hat{\gamma} - \gamma_0, \Pi_n v^* \rangle + \frac{1}{2}\|\epsilon_n \Pi_n v^*\|^2 + o_p(\epsilon_n n^{-1/2}) + o_p(n^{-1})$$
$$= \pm \epsilon_n \times \langle \hat{\gamma} - \gamma_0, v^* \rangle + \epsilon_n \times o_p(n^{-1/2}) + o_p(n^{-1}).$$

In summary, Claims 1, 2 and 3 imply that

$$0 \leq \frac{1}{n}\sum_{t=2}^{n}[\ell(\hat{\gamma}, Z_t) - \ell(\hat{\gamma} \pm \epsilon_n \Pi_n v^*, Z_t)]$$
$$= \mp \epsilon_n \frac{1}{n}\sum_{t=2}^{n}\frac{\partial \ell(\gamma_0, Z_t)}{\partial \gamma'}[v^*] \pm \epsilon_n \times \langle \hat{\gamma} - \gamma_0, v^* \rangle + \epsilon_n \times o_p(n^{-1/2}) + o_p(n^{-1})$$
$$= \mp \epsilon_n \mu_n\left(\frac{\partial \ell(\gamma_0, Z_t)}{\partial \gamma'}[v^*]\right) \pm \epsilon_n \times \langle \hat{\gamma} - \gamma_0, v^* \rangle + \epsilon_n \times o_p(n^{-1/2}) + o_p(n^{-1}).$$

Thus we obtain

$$\sqrt{n}\langle \hat{\gamma} - \gamma_0, v^* \rangle = \sqrt{n}\mu_n\left(\frac{\partial \ell(\gamma_0, Z_t)}{\partial \gamma'}[v^*]\right) + o_p(1) \Rightarrow N(0, \|v^*\|^2),$$

where the asymptotic normality is guaranteed by Billingsley's (1961a) ergodic stationary martingale difference CLT, and the asymptotic variance being equal to $\|v^*\|^2 \equiv \|\frac{\partial \rho(\gamma_0)}{\partial \gamma'}\|^2$ is implied by Lemma 4.1(1) and the definition of the Fisher norm $\|\cdot\|$. □

PROOF OF THEOREM 4.2. Given our normality results in Theorem 4.1, for our model we can take $\Sigma_n(v) = \frac{1}{\sqrt{n}}\sum_{t=2}^{n}\frac{\partial l(\gamma_0, Z_t)}{\partial \gamma'}[v]$, which is linear in $v$ and converges in distribution to $N(0, \|v\|^2)$, and $\frac{1}{2n}\sum_{t=2}^{n}(\frac{\partial l(\gamma_0, Z_t)}{\partial \gamma'}[v])^2 = \frac{1}{2}\|v\|^2 + o_p(1)$, and hence LAN holds. Notice that the proof in Wong (1992) allows for time series data, and following his proof, under LAN, we obtain that $\rho(\hat{\gamma}_n)$ achieves the semiparametric efficiency bound. Alternatively, we can conclude that $\rho(\hat{\gamma}_n)$ is semiparametrically efficient by applying the result of Bickel and Kwon (2001) which allows for strictly stationary semiparametric Markov models. □

PROOF OF PROPOSITION 4.1. Thanks to Lemma 4.1, we can directly extend the results in Bickel et al. (1993) for bivariate copula models with i.i.d. data to our copula-based first-order Markov time series setting. So the semiparametric efficiency bound for $\alpha_0$ is $\mathcal{I}_*(\alpha_0) = E_0\{\mathcal{S}_{\alpha_0}\mathcal{S}'_{\alpha_0}\}$, where $\mathcal{S}_{\alpha_0}$ is the *efficient score function* for $\alpha_0$, which is defined as the ordinary score function for $\alpha_0$ minus its population least squares orthogonal projection onto the closed linear span (clsp) of the score functions for the nuisance parameters $g_0$. And $\alpha_0$ is $\sqrt{n}$-*efficiently estimable* if and only if $E_0\{\mathcal{S}_{\alpha_0}\mathcal{S}'_{\alpha_0}\}$ is *nonsingular* [see e.g. Bickel et al. (1993)]. Hence (4.3) is clearly a necessary condition for $\sqrt{n}$-normality and efficiency of $\widehat{\alpha}_n$ for $\alpha_0$. Under Assumptions 4.2, 4.3 and 4.4', Propositions 4.7.4 and 4.7.6 of Bickel et al. [(1993), pages 165–168] for bivariate copula models apply. Therefore, with $\mathcal{S}_{\alpha_0}$ defined in (4.4), we have that $\mathcal{I}_*(\alpha_0) = E_0\{\mathcal{S}_{\alpha_0}\mathcal{S}'_{\alpha_0}\}$ is finite, positive-definite. This implies that Assumption 4.4 is satisfied with $\rho(\gamma) = \lambda'\alpha$ and $\omega = \infty$ and $\|v^*\|^2 = \|\frac{\partial \rho(\gamma_0)}{\partial \gamma'}\|^2 = \lambda'\mathcal{I}_*(\alpha_0)^{-1}\lambda < \infty$. By Theorem 4.1, for any $\lambda \in \mathcal{R}^d, \lambda \neq 0$, we have $\sqrt{n}(\lambda'\widehat{\alpha}_n - \lambda'\alpha_0) \Rightarrow \mathcal{N}(0, \lambda'\mathcal{I}_*(\alpha_0)^{-1}\lambda)$. This implies Proposition 4.1. □

PROOF OF THEOREM 5.1. The proof basically follows from that of Shen and Shi (2005), except for our definition of joint log-likelihood, our definition of Fisher norm $\|\cdot\|$, and our application of Billingsley's CLT for ergodic stationary martingale difference processes. These modifications are the same as those in our proof of Theorem 4.1. A detailed proof is omitted due to the length of the paper but is available upon request. □

## APPENDIX B: TABLES AND FIGURES

**Different estimators**: Sieve = Sieve MLE; Ideal = Ideal MLE; 2step = Chen–Fan; Para = correctly specified parametric MLE; Mis-N = parametric MLE using misspecified normal distribution as marginal; Mis-EV = parametric MLE using misspecified extreme value distribution as marginal.

Results are all based on 1000 MC replications of estimates using $n = 1000$ time series simulation. $\text{Bias}^2_{10^3}$, $\text{Var}_{10^3}$ and $\text{MSE}_{10^3}$ are the true values of Bias$^2$, Var and MSE multiplied by 1000 respectively. $\tau =$ Kendall's $\tau$, $\lambda =$ lower tail dependence index.

TABLE 1
*Clayton copula, true marginal $G = t_3$: estimation of $\alpha$*

|  |  | **Sieve** | **Ideal** | **2step** | **Para** | **Mis-N** | **Mis-EV** |
|---|---|---|---|---|---|---|---|
| $\alpha = 2$ | Mean | 1.969 | 2.002 | 1.912 | 1.989 | 2.400 | 2.957 |
| $\tau(0.500)$ | Bias | $-0.031$ | 0.002 | $-0.088$ | $-0.011$ | 0.400 | 0.957 |
| $\lambda(0.707)$ | Var | 0.019 | 0.007 | 0.101 | 0.012 | 0.103 | 0.056 |
|  | MSE | 0.020 | 0.007 | 0.109 | 0.012 | 0.264 | 0.971 |
|  | $\alpha^{\text{MC}}_{(2.5,97.5)}$ | (1.70, 2.25) | (1.83, 2.17) | (1.36, 2.60) | (1.76, 2.19) | (1.99, 3.28) | (2.57, 3.36) |

Table 1
*(Continued)*

|  |  | Sieve | Ideal | 2step | Para | Mis-N | Mis-EV |
|---|---|---|---|---|---|---|---|
| $\alpha = 5$ | Mean | 4.849 | 5.003 | 4.359 | 4.979 | 5.859 | 5.923 |
| $\tau(0.714)$ | Bias | $-0.151$ | 0.003 | $-0.642$ | $-0.021$ | 0.859 | 0.923 |
| $\lambda(0.871)$ | Var | 0.093 | 0.026 | 1.247 | 0.041 | 0.189 | 0.338 |
|  | MSE | 0.116 | 0.026 | 1.658 | 0.042 | 0.927 | 1.190 |
|  | $\alpha^{MC}_{(2.5,97.5)}$ | (4.25, 5.48) | (4.69, 5.32) | (2.67, 7.12) | (4.58, 5.35) | (5.36, 6.95) | (4.89, 6.62) |
| $\alpha = 10$ | Mean | 9.687 | 10.00 | 7.115 | 9.967 | 11.42 | 11.57 |
| $\tau(0.833)$ | Bias | $-0.313$ | 0.004 | $-2.886$ | $-0.033$ | 1.425 | 1.570 |
| $\lambda(0.933)$ | Var | 0.351 | 0.085 | 4.852 | 0.129 | 0.577 | 1.194 |
|  | MSE | 0.449 | 0.085 | 13.18 | 0.130 | 2.607 | 3.659 |
|  | $\alpha^{MC}_{(2.5,97.5)}$ | (8.68, 10.87) | (9.43, 10.6) | (3.87, 12.5) | (9.26, 10.6) | (10.33, 12.9) | (9.68, 12.9) |
| $\alpha = 12$ | Mean | 11.62 | 12.01 | 7.896 | 11.98 | 13.67 | 13.82 |
| $\tau(0.857)$ | Bias | $-0.382$ | 0.012 | $-4.104$ | $-0.016$ | 1.668 | 1.816 |
| $\lambda(0.944)$ | Var | 0.541 | 0.119 | 5.656 | 0.222 | 0.770 | 1.917 |
|  | MSE | 0.687 | 0.120 | 22.50 | 0.222 | 3.552 | 5.214 |
|  | $\alpha^{MC}_{(2.5,97.5)}$ | (10.5, 13.3) | (11.3, 12.7) | (4.35, 13.6) | (11.0, 12.9) | (12.3, 15.7) | (11.4, 15.4) |

**Acknowledgments.** We are grateful to the co-editor, an associate editor and two anonymous referees for their very helpful comments. We thank Don Andrews, Brendan Beare, Rohit Deo, Rustam Ibragimov, Per Mykland, Andrew Patton, Peter Phillips, Demian Pouzo, Zhijie Xiao for very useful discussions, and Zhentao Shi for excellent research assistance on Monte Carlo studies. We also acknowledge helpful comments from the participants at 2007 Exploratory Seminar on Stochastics and Dependence in Finance, Risk Management and Insurance at the Radcliffe Institute for Advanced Study at Harvard University, 2008 International Symposium on Financial Engineering and Risk Management (FERM), China, 2008 Conference in honor of Peter C. B. Phillips, SMU, 2008 Oxford-Man Institute Symposium on Modelling Multivariate Dependence and Extremes in Finance, and 2008 econometrics seminar at Columbia University.

TABLE 2
*Gumbel copula, true marginal $G = t_3$: estimation of $\alpha$*

|  |  | Sieve | Ideal | 2step | Para | Mis-N | Mis-EV |
|---|---|---|---|---|---|---|---|
| $\alpha = 2$ | Mean | 2.002 | 1.999 | 1.982 | 1.992 | 2.377 | 1.864 |
| $\tau(0.5)$ | Bias | 0.002 | $-0.001$ | $-0.018$ | $-0.008$ | 0.377 | $-0.136$ |
|  | Var | 0.007 | 0.002 | 0.013 | 0.005 | 0.153 | 0.026 |
|  | MSE | 0.007 | 0.002 | 0.014 | 0.005 | 0.295 | 0.045 |
|  | $\alpha^{\text{MC}}_{(2.5,97.5)}$ | (1.85, 2.18) | (1.91, 2.10) | (1.78, 2.23) | (1.85, 2.14) | (1.99, 3.55) | (1.60, 2.22) |
| $\alpha = 3.5$ | Mean | 3.486 | 3.498 | 3.352 | 3.481 | 3.906 | 3.629 |
| $\tau(0.714)$ | Bias | $-0.014$ | $-0.002$ | $-0.148$ | $-0.019$ | 0.406 | 0.129 |
|  | Var | 0.064 | 0.008 | 0.130 | 0.021 | 0.269 | 0.315 |
|  | MSE | 0.064 | 0.008 | 0.152 | 0.021 | 0.434 | 0.331 |
|  | $\alpha^{\text{MC}}_{(2.5,97.5)}$ | (3.06, 4.07) | (3.34, 3.68) | (2.76, 4.20) | (3.21, 3.87) | (3.21, 5.38) | (2.73, 4.83) |
| $\alpha = 6$ | Mean | 5.797 | 5.998 | 5.253 | 5.971 | 6.359 | 6.8805 |
| $\tau(0.833)$ | Bias | $-0.203$ | $-0.002$ | $-0.747$ | $-0.029$ | 0.359 | 0.881 |
|  | Var | 0.320 | 0.023 | 0.676 | 0.071 | 0.396 | 2.328 |
|  | MSE | 0.362 | 0.023 | 1.235 | 0.072 | 0.525 | 3.103 |
|  | $\alpha^{\text{MC}}_{(2.5,97.5)}$ | (4.67, 6.95) | (5.72, 6.31) | (3.92, 7.17) | (5.47, 6.67) | (5.20, 7.48) | (4.32, 9.78) |
| $\alpha = 7$ | Mean | 6.667 | 6.997 | 5.873 | 6.971 | 7.357 | 8.257 |
| $\tau(0.857)$ | Bias | $-0.333$ | $-0.003$ | $-1.127$ | $-0.029$ | 0.357 | 1.257 |
|  | Var | 0.456 | 0.032 | 0.968 | 0.106 | 0.506 | 3.859 |
|  | MSE | 0.566 | 0.032 | 2.238 | 0.107 | 0.633 | 5.438 |
|  | $\alpha^{\text{MC}}_{(2.5,97.5)}$ | (5.34, 8.12) | (6.67, 7.37) | (4.23, 8.20) | (6.34, 7.79) | (6.01, 8.58) | (4.96, 12.24) |

TABLE 3
*Clayton copula, true marginal $G = t_3$: estimation of G. Reported $Bias^2$, Var and MSE are the true ones multiplied by 1000*

|  |  | **Sieve** | | **2step** | | **Para** | | **Mis-N** | | **Mis-EV** | |
|---|---|---|---|---|---|---|---|---|---|---|---|
|  |  | $Q_{1/3}$ | $Q_{2/3}$ | $Q_{1/3}$ | $Q_{2/3}$ | $Q_{1/3}$ | $Q_{2/3}$ | $Q_{1/3}$ | $Q_{2/3}$ | $Q_{1/3}$ | $Q_{2/3}$ |
| $\alpha = 2$ | Mean | 0.325 | 0.673 | 0.333 | 0.666 | 0.333 | 0.667 | 0.347 | 0.557 | 0.382 | 0.614 |
| $\tau(0.500)$ | $Bias^2_{10^3}$ | 0.026 | 0.007 | 0.011 | 0.013 | 0.009 | 0.009 | 0.282 | 12.84 | 2.710 | 3.145 |
| $\lambda(0.707)$ | $Var_{10^3}$ | 0.054 | 0.049 | 1.430 | 0.801 | 0.002 | 0.002 | 1.921 | 5.651 | 0.755 | 0.947 |
|  | $MSE_{10^3}$ | 0.080 | 0.056 | 1.441 | 0.814 | 0.011 | 0.011 | 2.203 | 18.49 | 3.465 | 4.092 |
| $\alpha = 5$ | Mean | 0.322 | 0.671 | 0.332 | 0.667 | 0.333 | 0.667 | 0.331 | 0.537 | 0.342 | 0.579 |
| $\tau(0.714)$ | $Bias^2_{10^3}$ | 0.072 | 0.002 | 0.003 | 0.011 | 0.009 | 0.009 | 0.001 | 17.65 | 0.134 | 8.276 |
| $\lambda(0.871)$ | $Var_{10^3}$ | 0.081 | 0.085 | 6.474 | 2.969 | 0.002 | 0.002 | 1.401 | 5.697 | 2.234 | 5.346 |
|  | $MSE_{10^3}$ | 0.153 | 0.087 | 6.478 | 2.980 | 0.011 | 0.011 | 1.403 | 23.35 | 2.369 | 13.62 |
| $\alpha = 10$ | Mean | 0.319 | 0.664 | 0.331 | 0.666 | 0.333 | 0.667 | 0.364 | 0.584 | 0.371 | 0.624 |
| $\tau(0.833)$ | $Bias^2_{10^3}$ | 0.128 | 0.042 | 0.001 | 0.013 | 0.009 | 0.009 | 1.132 | 7.452 | 1.642 | 2.123 |
| $\lambda(0.933)$ | $Var_{10^3}$ | 0.109 | 0.137 | 22.28 | 9.800 | 0.003 | 0.003 | 0.711 | 3.410 | 2.103 | 4.192 |
|  | $MSE_{10^3}$ | 0.236 | 0.178 | 22.29 | 9.813 | 0.012 | 0.012 | 1.843 | 10.86 | 3.744 | 6.315 |
| $\alpha = 12$ | Mean | 0.318 | 0.661 | 0.331 | 0.665 | 0.333 | 0.667 | 0.374 | 0.598 | 0.375 | 0.633 |
| $\tau(0.857)$ | $Bias^2_{10^3}$ | 0.154 | 0.079 | 0.001 | 0.023 | 0.010 | 0.010 | 1.903 | 5.242 | 2.052 | 1.351 |
| $\lambda(0.944)$ | $Var_{10^3}$ | 0.127 | 0.141 | 28.83 | 12.08 | 0.003 | 0.003 | 0.950 | 2.662 | 2.494 | 4.934 |
|  | $MSE_{10^3}$ | 0.281 | 0.220 | 28.83 | 12.10 | 0.013 | 0.013 | 2.853 | 7.904 | 4.547 | 6.286 |

TABLE 4
*Gumbel copula, true marginal $G = t_3$: estimation of G*

|  |  | **Sieve** | | **2step** | | **Para** | | **Mis-N** | | **Mis-EV** | |
|---|---|---|---|---|---|---|---|---|---|---|---|
|  |  | $Q_{1/3}$ | $Q_{2/3}$ | $Q_{1/3}$ | $Q_{2/3}$ | $Q_{1/3}$ | $Q_{2/3}$ | $Q_{1/3}$ | $Q_{2/3}$ | $Q_{1/3}$ | $Q_{2/3}$ |
| $\alpha = 2$ | Mean | 0.328 | 0.673 | 0.333 | 0.666 | 0.333 | 0.667 | 0.401 | 0.613 | 0.519 | 0.737 |
| $\tau(0.500)$ | $Bias^2_{10^3}$ | 0.004 | 0.011 | 0.007 | 0.018 | 0.009 | 0.009 | 5.069 | 3.239 | 35.53 | 4.456 |
|  | $Var_{10^3}$ | 0.059 | 0.063 | 0.755 | 1.025 | 0.003 | 0.003 | 2.389 | 3.111 | 10.44 | 7.202 |
|  | $MSE_{10^3}$ | 0.063 | 0.074 | 0.762 | 1.043 | 0.012 | 0.012 | 7.457 | 6.350 | 45.98 | 11.66 |
| $\alpha = 3.5$ | Mean | 0.328 | 0.675 | 0.332 | 0.665 | 0.333 | 0.667 | 0.524 | 0.719 | 0.565 | 0.746 |
| $\tau(0.714)$ | $Bias^2_{10^3}$ | 0.004 | 0.025 | 0.005 | 0.030 | 0.009 | 0.009 | 37.55 | 2.386 | 55.42 | 5.762 |
|  | $Var_{10^3}$ | 0.139 | 0.147 | 2.353 | 3.482 | 0.004 | 0.004 | 18.71 | 9.238 | 28.40 | 18.12 |
|  | $MSE_{10^3}$ | 0.143 | 0.171 | 2.358 | 3.511 | 0.013 | 0.013 | 56.26 | 11.62 | 83.82 | 23.88 |
| $\alpha = 6$ | Mean | 0.325 | 0.681 | 0.330 | 0.664 | 0.333 | 0.667 | 0.501 | 0.700 | 0.497 | 0.676 |
| $\tau(0.833)$ | $Bias^2_{10^3}$ | 0.025 | 0.120 | 0.000 | 0.036 | 0.009 | 0.009 | 29.17 | 0.899 | 27.97 | 0.037 |
|  | $Var_{10^3}$ | 0.273 | 0.255 | 6.840 | 10.37 | 0.005 | 0.005 | 40.49 | 20.60 | 40.98 | 29.81 |
|  | $MSE_{10^3}$ | 0.298 | 0.375 | 6.840 | 10.41 | 0.014 | 0.014 | 69.66 | 21.50 | 68.96 | 29.84 |
| $\alpha = 7$ | Mean | 0.324 | 0.684 | 0.329 | 0.665 | 0.333 | 0.667 | 0.477 | 0.679 | 0.476 | 0.655 |
| $\tau(0.857)$ | $Bias^2_{10^3}$ | 0.041 | 0.182 | 0.000 | 0.029 | 0.009 | 0.009 | 21.46 | 0.076 | 21.35 | 0.227 |
|  | $Var_{10^3}$ | 0.314 | 0.275 | 9.362 | 13.79 | 0.006 | 0.006 | 49.51 | 26.89 | 45.82 | 33.93 |
|  | $MSE_{10^3}$ | 0.355 | 0.457 | 9.362 | 13.82 | 0.016 | 0.016 | 70.97 | 26.96 | 67.16 | 34.16 |

TABLE 5
*Clayton copula, true marginal $G = t_3$: estimation of 0.01 conditional quantile*

|  |  | **Sieve** | **Ideal** | **2step** | **Para** | **Mis-N** | **Mis-EV** |
|---|---|---|---|---|---|---|---|
| $\alpha = 5$ | IntBias$^2_{10^3}$ | 36.26 | 0.000 | 150.0 | 0.172 | 900.7 | 704.8 |
| $\tau(0.714)$ | IntVar$_{10^3}$ | 32.15 | 5.450 | 985.3 | 10.18 | 463.7 | 313.4 |
| $\lambda(0.871)$ | IntMSE$_{10^3}$ | 68.41 | 5.450 | 1135 | 10.35 | 1364 | 1018 |
| $\alpha = 10$ | IntBias$^2_{10^3}$ | 7.712 | 0.000 | 527.3 | 0.040 | 815.3 | 427.4 |
| $\tau(0.833)$ | IntVar$_{10^3}$ | 19.36 | 2.475 | 855.3 | 3.716 | 361.7 | 202.7 |
| $\lambda(0.933)$ | IntMSE$_{10^3}$ | 27.07 | 2.475 | 1383 | 3.756 | 1177 | 630.1 |
| $\alpha = 12$ | IntBias$^2_{10^3}$ | 2.851 | 0.000 | 367.7 | 0.004 | 181.1 | 175.9 |
| $\tau(0.857)$ | IntVar$_{10^3}$ | 6.236 | 1.068 | 590.9 | 1.578 | 59.44 | 46.12 |
| $\lambda(0.944)$ | IntMSE$_{10^3}$ | 9.086 | 1.069 | 958.7 | 1.582 | 240.5 | 222.0 |

For each $\alpha$, evaluation is based on the common support of 1000 MC simulated data. Reported integrated Bias$^2$, integrated Var and the integrated MSE are the true ones multiplied by 1000.

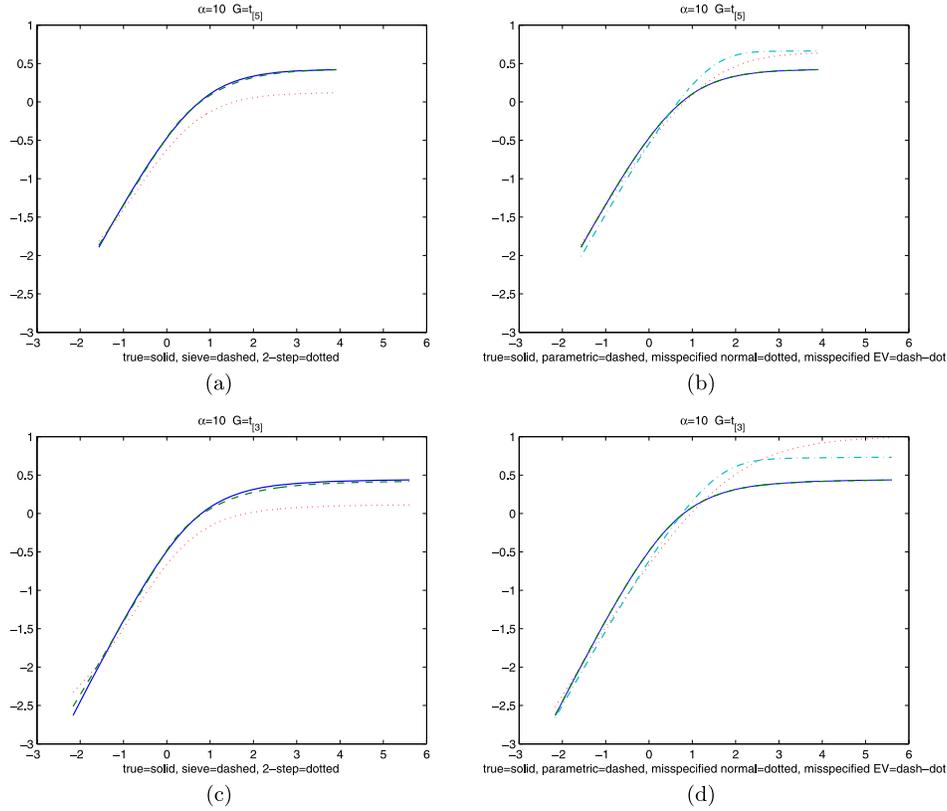

FIG. 2. *Clayton copula (true $\alpha = 10$, marginal $G = t_5$, $t_3$): estimation of 0.01 conditional quantile. Evaluation is based on the common support of 1000 MC simulated data.*

X. Chen
Cowles Foundation for Research in Economics
Yale University
30 Hillhouse Ave., Box 208281
New Haven, Connecticut 06520
USA
E-mail: xiaohong.chen@yale.edu

W. B. Wu
Department of Statistics
University of Chicago
5734 S. University Ave.
Chicago, Illinois 60637
USA
E-mail: wbwu@galton.uchicago.edu

Y. Yi
Department of Economics
New York University
19 West 4th St.
New York, New York 10012
USA
E-mail: yanping.yi@nyu.edu